\documentclass[11pt,a4paper]{article}
\usepackage[utf8]{inputenc}
\usepackage{amsmath,bm}
\usepackage{amsmath}
\usepackage{amsthm}
\usepackage{amsfonts}
\usepackage[square, comma, sort&compress,numbers]{natbib}
\usepackage{amssymb}
\usepackage{setspace}
\usepackage{graphicx}
\usepackage{color}
\usepackage[marginal]{footmisc}

\usepackage{indentfirst}
\usepackage[left=2cm,right=2cm,top=2cm,bottom=2cm]{geometry}
\author{}
\usepackage{float}
\newtheorem{cor}{Corollary}
\newtheorem{lem}{Lemma}
\newtheorem{thm}{Theorem}

\begin{document}
\bibliographystyle{plainnat}

\title{\textbf {Existence of Heterodimensional Cycles near Shilnikov Loops in Systems with a $\mathbb{Z}_2$ Symmetry}}
\date{}
\maketitle
\vspace{-1.5cm}
\begin{center}
{\large 
DONGCHEN LI \\[5pt]
{\it Department of Mathematics, Imperial College London\\
180 Queen's Gate, London SW7 2AZ, United Kingdom}\\[20pt]
DMITRY V. TURAEV\\[5pt]
{\it Department of Mathematics, Imperial College London\\
180 Queen's Gate, London SW7 2AZ, United Kingdom}\\[5pt]
{\it Department of Mathematics, Lobachevsky State
University of Nizhny Novgorod\\
23 Prospekt Gagarina, Nizhny Novgorod 603950, Russia}\\[5pt]
{\it Joseph Meyerhoff Visiting Professor, Weizmann Institute of Science\\ 234 Herzl Street, Rehovot 7610001, Israel}
}
\end{center}
\vspace{1cm}
\footnote{\noindent Date: October 16, 2016.\\
\noindent This work was supported by grant RSF 14-41-00044. The authors also acknowledge support by the Royal Society grant IE141468 and EU Marie-Curie IRSES Brazilian-European partnership in Dynamical Systems (FP7-PEOPLE-2012-IRSES 318999 BREUDS).
}
\par{}
\noindent{\bf Abstract.} We prove that a pair of heterodimensional cycles can be born at the bifurcations of a pair of Shilnikov loops (homoclinic loops to a saddle-focus equilibrium) having a one-dimensional unstable manifold in a volume-hyperbolic flow with a $\mathbb{Z}_2$ symmetry. We also show that these heterodimensional cycles can belong to a chain-transitive attractor of the system along with persistent homoclinic tangency.
\par{}
\par{}
\section{Introduction}
There is a point of view that the main feature of the dynamics of non-hyperbolic chaotic systems is the persistent coexistence of orbits with different numbers of positive Lyapunov exponents. This can be caused by the existence of either a homoclinic tangency or a heterodimensional cycle, i.e. a cycle which includes heteroclinic connections between saddle periodic orbits with different indices (dimensions of their unstable manifolds). In our opinion, for multidimensional systems (i.e. diffeomorphisms with dimension three or higher and flows with dimension four or higher), the most basic mechanism of this phenomenon must be a heterodimensional cycle. Heterodimensional cycles of co-index 1 were first studied by Newhouse and Palis in \citep{np76}. Here co-index is the difference between the indices of the corresponding periodic orbits of a heterodimensional cycle.
The fact that the non-transverse heteroclinic intersections in such cycles can be persistent was discovered by Díaz and collaborators (see \citep{d92,d95,d95_2,bd96});
a comprehensive theory of $C^1$-generic properties of diffeomorphisms having heterodimensional cycles of co-index 1 was built mostly in the works
of Bonatti and Díaz (see \citep{bd96,bd08}).
\par{}
In this paper we consider heterodimensional cycles for $C^r$ flows in $\mathbb{R}^n$ $(r\geqslant 3,n\geqslant 4)$. We give an example of a simple codimension-one homoclinic bifurcation which (among other things) results in the 
emergence of heterodimensional cycles of co-index 1. This is a symmetric version of the Shilnikov bifurcation of two homoclinic
loops to a saddle-focus. The symmetry also links heterodimensional cycles to Lorenz-like systems. We show in Section \ref{sec:1.2} an example of the system which satisfies the assumptions in this paper by adding an extra direction to the geometric Lorenz model. Namely, the system 
\begin{equation*}
\left\{
\begin{array}{rcl}
\dot{x} &=& \sigma (y-x),\\
\dot{y} &=& x(r-z) -y,\\
\dot{z} &=& -bz + xy + \varepsilon u, \\
\dot{u} &=& -(b+f)u - \varepsilon z,
\end{array}\right.
\end{equation*} 
\noindent will, for some choice of parameter values and function $f$, undergo the bifurcation which give rise to heterodimensional cycles (see Section \ref{sec:1.2} for more details ).
\par{}
It has been shown in \citep{os87} that under certain (open) condition on the eigenvalues of
the saddle-focus equilibrium the bifurcation of a homoclinic loop to the saddle-focus in three-dimensional systems creates coexisting saddles periodic orbits with different indices. We generalise this result for systems with dimension four or higher (see Corollary \ref{cor:dense} in Section \ref{sec:prf1.3}) and show that a symmetric pair of such loops can be split in such a way that some of these saddles
acquire heteroclinic connections and the heterodimensional cycles are formed (see Theorem \ref{thm:hetero_1} in Section \ref{sec:2.1}). One should note that it is impossible to create a heterodimensional cycle only using saddles near one single homoclinic loop (under condition C1 in Section \ref{sec:1.1}). The interplay of two homoclinic loops is crucial and we will explain this in Section \ref{sec:1.1}.
\par{}
We remark here that, by imposing the symmetry requirement, the codimension of the bifurcation under consideration is brought down to one. This is because that the existence of one homoclinic loop now implies the existence of the second one; moreover, the coincidence condition (which is an equality-type condition) needed for the emergence of heterodimensional cycles will be fulfilled automatically (see Section \ref{sec:1.1}). In other words, the symmetry allows us to give a relatively simple criterion for the heterodimensional cycle chaos. More specifically, under the symmetry condition, the appearance of a single Shilnikov loop with the volume-hyperbolicity near the equilibrium (condition C3 in Section \ref{sec:1.1}) is sufficient to show the existence of heterodimensional cycles in systems which can be arbitrarily close to the original one. Therefore, we can obtain a complex structure from a simple one. The computations we do here are quite involved. This is caused by the fact that we need to consider the perturbations which keep the symmetry of the system. Creating heterodimensional cycles for general systems without the symmetry is easier, but the bifurcations become codimension three; this case is considered in \citep{li16}.
\par{}
We also show that if the Shilnikov loops are originally within an attractor, then the heterodimensional cycles obtained in this paper can belong to this attractor (see Theorem \ref{thm:hetero_2} in Section \ref{sec:2.2}). The attractor considered here is the one proposed in \citep{ts98} which is chain-transitive volume-hyperbolic, and contains the equilibrium and its two separatrices.
\subsection{Problem Setting}\label{sec:1.1}
\par{}
In what follows we describe the system considered in this paper, and give conditions required to create heterodimensional cycles via homoclinic bifurcations. 
\par{}
Let us consider a $C^r$ flow $X$ in $\mathbb{R}^n$ (where $r\geqslant 3, n \geqslant 4$) having an equilibrium $O$ with a one-dimensional unstable manifold and a homoclinic loop associated to $O$. We assume that system $X$ satisfies the conditions below.
\\[10pt]
\indent (C1) (Non-degeneracy condition) The extended unstable manifold $W^{uE}(O)$ 
is transverse to the strong-stable foliation $\mathcal{F}_0$ of the stable manifold $W^s(O)$ at the points of the homoclinic loop.
\\[10pt]
\indent An extended unstable manifold $W^{uE}(O)$  is a smooth three-dimensional invariant manifold which is tangent at the points of $W^u_{loc}$ to the eigenspace corresponding to the unstable and weak stable characteristic exponents (those closest to the imaginary axis from right). It contains the stable manifold $W^u(O)$ 
and is transverse to the strong-stable manifold $W^{ss}_{loc}(O)$ at $O$. The foliation $\mathcal{F}_0$ is the uniquely defined, smooth, invariant foliation of the stable manifold, which includes $W^{ss}(O)$ as one of its leaves. We will discuss more on this foliation later.
\par{}
Note that condition (C1) is open and dense in $C^r$ topology, i.e., 
if it is not fulfilled initially, then it can be achieved by an arbitrarily small perturbation of the system; once this condition is satisfied, it holds for every $C^r$-close system. We proceed to listing other conditions.
\\[10pt]
\indent (C2) The equilibrium $O$ is a saddle-focus, and the eigenvalues of the linearised matrix of $X$ at $O$ are $
\gamma,-\lambda+\omega i,-\lambda-\omega i,\alpha_j $
such that
\begin{equation*}
\omega\neq 0 \quad \mbox{and}\quad \mbox{Re }(\alpha_j)<-\lambda<0<\gamma \,\, (j=1,2\ldots n-3).
\end{equation*}
\indent It follows from the result in Appendix A of \citep{sstc1} that if (C2) is satisfied, then system $X$ near $O$ can be brought to the form
\begin{equation}\label{eq:1setting_1}
\begin{array}{rcl}
\dot{x}&=&x, \\[5pt]
\dot{y_1}&=&-\rho y_1 - \omega y_2 + f_{11}(x,y,z)y+f_{12}(x,y,z)z, \\[5pt]
\dot{y_2}&=&\omega y_1 - \rho y_2 + f_{21}(x,y,z)y+f_{22}(x,y,z)z, \\[5pt]
\dot{z}  &=&Bz+f_{31}(x,y,z)y+f_{32}(x,y,z)z,    
\end{array} 
\end{equation}
\noindent by some $C^{r-1}$-transformation of coordinates and time (with assuming $\gamma=1$). Here $x=(x_1,x_2)$, and the eigenvalues of matrix $B$ are $\alpha_1\ldots \alpha_{n-3}$.
Functions $f_{ij}$ are $C^{r-1}$ smooth and satisfy
\begin{equation}\label{eq:1setting_2}
f_{ij}(0,0,0)=0,\,\, f_{1j}(0,y,z)\equiv 0,\,\, f_{2j}(0,y,z)\equiv 0, \,\, f_{i1}(x,0,0) \equiv 0 \,\,(i=1,2,3;j=1,2).
\end{equation} 
In such coordinate system, the coordinates of $O$ are $(0,0,0)$ and the local invariant manifolds are straightened, i.e. we have 
\begin{equation*}
W^u_{loc}(O)=\{y=0,z=0\},\quad W^s_{loc}(O)=\{x=0\},\quad W^{ss}_{loc}(O)=\{x=0,y=0\}.
\end{equation*}
\par{}
The one-dimensional unstable manifold of $O$ consists of two separatrices; the upper one, $\Gamma^+$ corresponds, locally, to $x>0$
and the lower separatrix $\Gamma^-$ corresponds to $x<0$. Let the upper separatrix return to $O$ as $t\to+\infty$ and form a homoclinic
loop. Thus, the homoclinic loop, when it tends to $O$ as $t=-\infty$, coincides with a piece of the $x$-axis, and when the loop tends to $O$ as $t\to+\infty$
it lies in $\{x=0\}$. In this coordinate system the extended unstable manifold $W^{uE}(O)$ is tangent to $\{z=0\}$ at points of $\Gamma^+$.
\par{}
We also impose a condition which ensures that the system expands three-dimensional volumes near $O$:
\\[10pt]
\indent (C3) The ratio $\rho=\dfrac{\lambda}{\gamma}<\dfrac{1}{2}$.
\par{}
Under condition (C1), conditions (C2) and (C3) are necessary for obtaining heterodimensional cycles via bifurcations of homoclinic loops. Indeed, if $O$ is a saddle or $O$ is a saddle-focus with $\rho>1$, then at most one periodic orbit can be born from the bifurcation of one homoclinic loop (see \citep{sstc2,hs10}); in the case where $O$ is a saddle-focus with $1/2<\rho<1$, there can be infinitely many coexisting periodic orbits with indices 1 and 2 near one homoclinic loop (see \citep{os87}), but index-1 orbits are attractors, i.e. not saddles, and therefore they cannot be used to create heterodimensional cycles.
\par
\begin{figure}[!h]
\begin{center}
\includegraphics[width=0.5\textwidth]{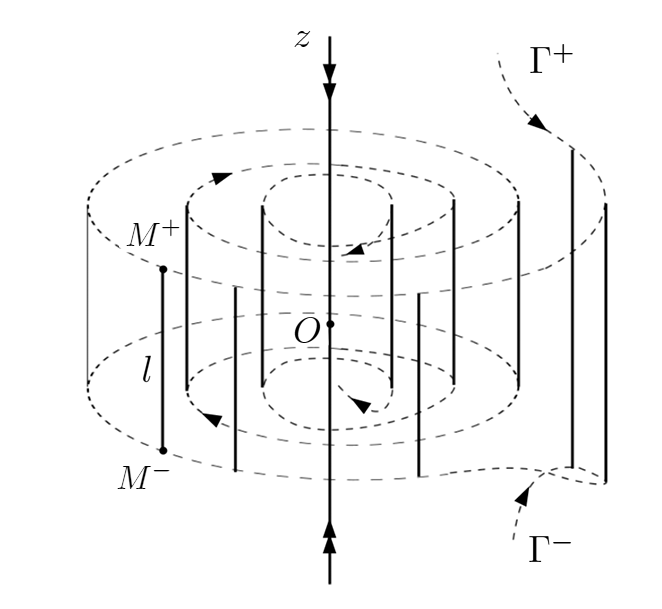}
\end{center}
\caption{The dashed curves represent the two homoclinic loops and the solid vertical lines represent leaves of the foliation $\mathcal{F}_0$. The coincidence condition for system $X$ is that, for any point of $\Gamma^+$ lying in a leaf $l$, there exists one point of $\Gamma^-$ that also lies in $l$. }
\label{fig:coincidence1}
\end{figure}

{}
When $O$ is a saddle-focus with $0<\rho<1/2$, one can obtain infinitely many coexisting saddle periodic orbits of indices 2 and 3 near one homoclinic loop (see Lemma \ref{lem:dense_i2} and Corollary \ref{cor:dense} generalize the same result obtained for three-dimensional systems in \citep{os87}). However,  it is known (see \citep{tu96}) that, under some genericity assumption on the homoclinic loop $\Gamma^+$ (i.e. condition C1 of this paper), the original system, and every system close to it, has a three-dimensional invariant manifold $\mathcal{M}$ such that every orbit which lies entirely in a small neighbourhood of $O\cap\Gamma^+$ must lie in $\mathcal{M}$. This gives a robust three-dimensional reduction of the dynamics near $\Gamma^+$, which prevents the birth of heterodimensional cycles at any bifurcations of $\Gamma^+$. Therefore, the interplay of two homoclinic loops is required, and we assume that the separatrix $\Gamma^-$ also forms a homoclinic loop. Indeed, the existence of essentially four-dimensional dynamics can be guaranteed if system $X$ satisfies the coincidence condition: the loops $\Gamma^+$ and $\Gamma^-$ intersect the same set of leaves of the strong-stable foliation $\mathcal{F}_0$ on $W^s(O)$ (see \citep{tu96,bc15}). This means that, for any point $M^+\in\Gamma^+$ lying in a leaf $l$, there exists a point $M^-\in\Gamma^-$ lying in the same leaf $l$ (see figure \ref{fig:coincidence1}).
\par{}
We now achieve this coincidence condition by imposing the symmetry:\\[10pt]
\indent (C4) System $X$ is invariant with respect to the transformation $R: (x,y,z)\to(-x,y,{\cal S} z)$
where ${\cal S}$ is a non-trivial involution which changes signs of some of the $z$-coordinates.
\par{}
With this condition satisfied, the existence of the loop $\Gamma^+$ implies the existence of the second homoclinic loop $\Gamma^-$, and moreover
the $y$-component will be the same for both of these homoclinic solutions. Besides, the above-mentioned coincidence is fulfilled automatically. In the rest of this paper we show that bifurcations of this pair of homoclinic loops can lead to the birth of heterodimensional cycles of co-index 1.
\subsection{An example}\label{sec:1.2}
A concrete example of a system satisfying conditions (C1) - (C4) can be found by a modification of the well-known Lorenz model given by 
\begin{equation*}
\left\{
\begin{array}{rcl}
\dot{x} &=& \sigma (y-x),\\
\dot{y} &=& x(r-z) -y,\\
\dot{z} &=& -bz + xy.
\end{array}\right.
\end{equation*}
\noindent Here the unstable and strong-stable directions are given by linear combinations of $x$ and $y$, and the weak-stable direction corresponds to coordinate $z$. This system is symmetric with respect to the transformation $(x,y,z)\to(-x,-y,z)$. We add a new variable $u$, and consider the system of the form
\begin{equation}\label{eq:example}
\left\{
\begin{array}{rcl}
\dot{x} &=& \sigma (y-x),\\
\dot{y} &=& x(r-z) -y,\\
\dot{z} &=& -bz + xy + \varepsilon u, \\
\dot{u} &=& -(b+f(x,y,z,u))u - \varepsilon z,
\end{array}\right.
\end{equation} 
\noindent where $f$ can be any non-linear function such that
\\[5pt]
\indent (1) the new system satisfies the symmetry with respect to the transformation $(x,y,z,u)\to(-x,-y,z,u)$ (which is the same as the $R$ symmetry introduced before); and
\\[5pt]
\indent (2) we have $f(0,0,0,0)=\partial f(0,0,0,0)/\partial (x,y,z,u)=0$, and the sum $(b+f)$ is close to zero outside a small neighbourhood of the equilibrium $(0,0,0,0)$.
\par{}
It is known (\citep{tu02,abs77,abs83}) that there is an open set in the parameter space around $(\sigma=10,b=8/3,r=28)$ such that, for parameter values inside this set, the Lorenz system has a strong-stable foliation and the two-dimensional areas transverse to the foliation are expanding near the equilibrium. The above property (2) ensures that the strong-stable foliation is inherited by the new system \eqref{eq:example}, and three-dimensional volumes transverse to the strong-stable foliation are expanded which implies condition (C3). Besides, property (2) also leads to the existence of an absorbing domain containing the equilibrium which implies the existence of a volume-hyperbolic attractor inside the domain (see Section \ref{sec:2.2} for more details). We note that the characteristic exponents corresponding to coordinates $z$ and $u$ are conjugate complex numbers, and therefore the original equilibrium in the Lorenz model now becomes a saddle-focus. There is numerical evidence (e.g. \citep{bss12}) that, for a dense subset of parameter values inside the open set near $(\sigma=10,b=8/3,r=28)$, the Lorenz system has a symmetric pair of homoclinic loops. The loops will persist after we add the extra coordinate $u$. Thus, at least for certain parameter values, this new system falls into the class of systems considered in this paper. 
\par{}
By Theorem \ref{thm:hetero_1} of this paper, we can have heterodimensional cycles by perturbing this system. More specifically, there exist certain parameter values of $\sigma,b$ and $r$, and function $f$ such that, for any $\delta>0$, one can find $C^r$ functions $g_1,g_2,g_3$ and $g_4$ with $\|g_i\|_{C^r}<\delta$ such that the system
\begin{equation*}
\left\{
\begin{array}{rcl}
\dot{x} &=& \sigma (y-x)+g_1,\\
\dot{y} &=& x(r-z) -y+g_2,\\
\dot{z} &=& -bz + xy + \varepsilon u+g_3, \\
\dot{u} &=& -(b+f(x,y,z,u))u - \varepsilon z +g_4,
\end{array}\right.
\end{equation*} 
has a symmetric pair of heterodimensional cycles. Moreover, these cycles belong to the above-mentioned attractor (see Theorem \ref{thm:hetero_2}).
\section{Results}\label{sec:2}
\subsection{Birth of heterodimensional cycles}\label{sec:2.1}
The main result of this paper is the following:
\begin{thm}\label{thm:hetero_1}
If system $X$ satisfies conditions (C2) - (C4), then in any arbitrarily small $C^r$ neighbourhood $(r\geqslant 1)$ of $X$ in the space of the $R$-symmetric systems, there exists a system which has
a symmetric pair of homoclinic loops to $O$, and a symmetric pair of heterodimensional cycles near these loops. Each heterodimensional cycle is associated to two periodic orbits of indices 2 and 3.
\end{thm}
\par{}
Condition (C1) is not mentioned here since it can also be obtained by an arbitrarily small perturbation (without destroying the loops). We remark here that we need $r\geqslant 3$ in our computations. If system $X$ is originally $C^r$ with $r=1$ or 2, then we can first make it $C^\infty$ by an arbitrarily small perturbation in $C^r$ topology, and recover the homoclinic loops (if destroyed) by an additional arbitrarily $C^r$-small perturbation. After this, we can perturb the system again to create heterodimensional cycles. We shall mention that, in the proof of the non-empty quasi-transverse intersection in Section \ref{sec:prf1.4}, we need to use the smooth dependence of the invariant manifolds of periodic orbits on the right-hand side of system \ref{eq:1setting_1}, and this is allowed after we make our system $C^\infty$.
\par{}
Before we sketch the proof, let us discuss more on the strong-stable foliation $\mathcal{F}_0$. In the coordinates of (\ref{eq:1setting_1}), the leaves of $\mathcal{F}_0$ on $W^s(O)_{loc}$ are given by $(x=0,y=\mbox{const})$. The non-degeneracy condition (C1) implies that the closed invariant set $O\cup\Gamma^+\cup\Gamma^-$ is partially hyperbolic: at the points of this set the contraction along the 
strong-stable leaves is stronger than a possible contraction in the directions tangent to $W^{uE}$. The partial hyperbolicity implies that the strong-stable foliation $\mathcal{F}_0$ extends (see \citep{hps77,tu96,ts98}), as a locally invariant, absolutely continuous foliation with smooth leaves, to a neighbourhood $U$
of $O\cup\Gamma^+\cup\Gamma^-$, and the foliation persists for all $C^r$-close systems. See \citep{an67,hps77} for more details on the properties of such foliation.
\par
\begin{figure}[!h]
\begin{center}
\includegraphics[width=0.5\textwidth]{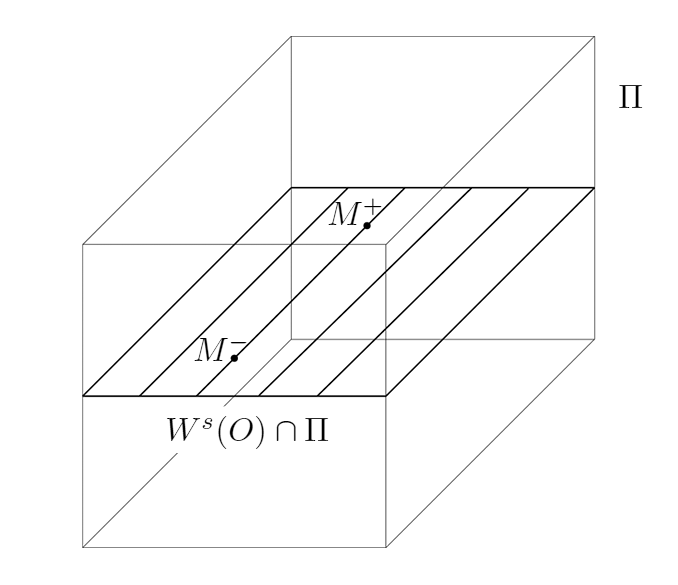}
\end{center}
\caption{For a four-dimensional system, the intersection points $M^+$ and $M^-$ on a small three-dimensional cross-section $\Pi$ belong to the same leaf of the foliation $\mathcal{F}_1$.}
\label{fig:coincidence2}
\end{figure}

{}
We take a small cross-section $\Pi$ transverse to the local stable manifold $W^s_{loc}(O)$ such that both loops $\Gamma^+$ and $\Gamma^-$ intersect $\Pi$. The flow induces a Poincaré map $T$ on $\Pi$. The intersections of the orbits of the leaves of $\mathcal{F}_0$ by the flow with the cross-section $\Pi$ form a strong-stable invariant foliation $\mathcal{F}_1$ for the Poincaré map $T$, which has leaves of the form $(x,y)=h(z)$ where the derivative $h'(z)$ is uniformly bounded. The detailed sufficient condition for the existence of such strong-stable foliation is proposed in \citep{ts98} and our system $X$ satisfies this condition. Note that the coincidence condition (given by the symmetry) implies that the projections of $\Gamma^+$ and $\Gamma^-$ onto any transversal along leaves of $\mathcal{F}_0$ coincide. Therefore, the intersection points $M^+$ and $M^-$ of $\Gamma^+$ and $\Gamma^-$ with $\Pi$ lie on the same leaf of $\mathcal{F}_1$, and have the same $y$-coordinate since $\Pi$ is near $O$ and the foliations on $W^s_{loc}(O)$ are straightened (see figure \ref{fig:coincidence2}). 
\par{}
The foliation $\mathcal{F}_1$ is invariant such that $T^{-1}(l\cap(T(\Pi)))$ is a leaf of the foliation if the intersection is non-empty. The foliation is also contracting in the sense that, for any two points in the same leaf, the distance between their iterates under the map $T$ tends to zero exponentially. Besides, by the absolute continuity of the foliation, the projection along the leaves from one transversal to another one changes areas by a finite multiple bounded away from zero. Note that the condition $\rho<{1}/{2}$
implies that the flow near $O$ expands three-dimensional volume in the $(x,y)$-space; the partial hyperbolicity of the flow near $O\cup\Gamma^+\cup\Gamma^-$ and the fact that the
orbits in $U$ spend only a finite time between successive returns to the small neighbourhood of $O$ imply that the flow in $U$ uniformly
expands the three-dimensional volume transverse to the strong-stable foliation (see \citep{tu96,ts98}). Consequently, the Poincaré map $T$ expands two-dimensional
areas transverse to the strong-stable foliation on $\Pi$. 
\par{}
In what follows, we describe the steps for creating a heterodimensional cycle in a small neighbourhood of $\Gamma^+\cap\Gamma^-\cap O$. At the birth of this cycle, a second one will be obtained automatically by the symmetry.
\par{}
First, we note that, according to Shilnikov theorem (see \citep{sh65,sh70}), each of the homoclinic loops $\Gamma^+$ and $\Gamma^-$ is accumulated by a countable set of single-round index-2 saddle periodic orbits of the flow (we call the orbit of the flow $n$-round if
it intersects $\Pi$ exactly $n$ times). Consequently, there exist two sets $\{P^+_k\}$ and $\{P^-_k\}$ of index-1 saddle fixed points of $T$ in $\Pi$ such that $P^+_k\to M^+$ and $P^-_k\to M^-$ as $k\to +\infty$. Any finite number of these points survive sufficiently small perturbations of the system. We embed $X$ into a two-parameter family $X_{\mu,\rho}$ of $R$-symmetric systems such that homoclinic loops split with a non-zero velocity as $\mu$ changes. As before, we take $\mu$ as the $x$-coordinate of the point $M^+$ where the upper separatrix $\Gamma^+$ first intersects $\Pi$ (so $-\mu$ is the $x$-coordinate
of the point $M^-$ of the first intersection of $\Gamma^-$ with $\Pi$). The second parameter is the ratio $\rho={\lambda}/{\gamma}$. It is well-known that arbitrarily close to $\mu=0$ there are values of $\mu$ for which both $\Gamma^+$ and
$\Gamma^-$ form a double-round homoclinic loop (see \citep{eff82,fe93,ga83,gtgn97}). Crucially, we show that by an arbitrarily small perturbation of $\rho$ (in addition to that of $\mu$), at the moment of existence of two double-round homoclinic loops, the unstable manifold of a point $P$ from the set $\{P_k^+\}$ that survives the splitting of the original loop intersects the strong-stable manifold of the point $M^-$ (see Lemma \ref{lem:pm}).
\par
\begin{figure}[!h]
\begin{center}
\includegraphics[width=1\textwidth]{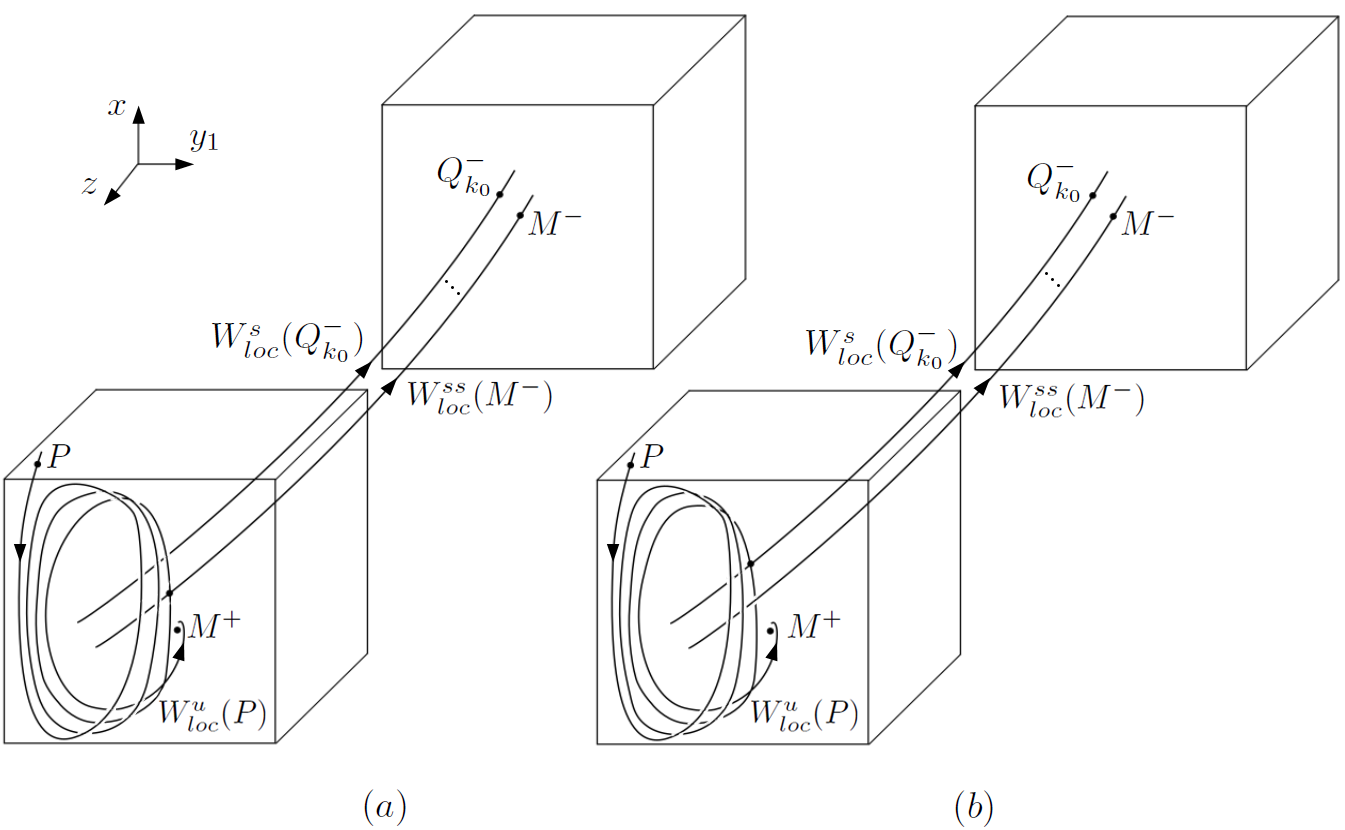}
\end{center}
\caption{As shown in figure (a), we can create an infinite sequence of index-2 point $Q_k^-$ accumulating on $M^-$ while keeping the intersection $W^u(P)\cap W^{ss}(M^-)$ by changing $\mu,\rho$ and $\nu$ together. In figure (b), the intersection $W^u(P)\cap W^s(Q_{k_0}^-)$ is created by changing $\nu$).}
\label{fig:3d_3}
\end{figure}

{}
Next, we use a generalization of Theorem 3 of \citep{os87} that if a system has a homoclinic loop to a saddle-focus with $\rho<1/2$,
then, by an arbitrarily perturbation which changes the value of $\rho$ without splitting the loop, one can create an infinite sequence of 
double-round saddle periodic orbits with three-dimensional unstable manifold which converges to the loop (see Lemma \ref{lem:dense_i2}). In our situation, we 
can consider a family of perturbations localized in a sufficiently small neighbourhood of $O$ such that neither the symmetry of the system is broken,
nor the double-round loops are split, nor the heteroclinic intersection between $W^u(P_k^+)$ and $W^{ss}(M^-)$ is destroyed, while
the value of $\rho$ changes with a non-zero velocity. Then, at an appropriately chosen value of $\rho$ the double-round loop $\Gamma^-$ becomes a limit of a sequence of 4-round saddle periodic orbits with three-dimensional unstable manifold (see Figure \ref{fig:3d_3} (a)).  On the cross-section $\Pi$, we thus have
an infinite sequence of index-2 saddle points $Q_k^-$ of period 4 which converges to $M^-$; the stable manifolds of these points are given by
the leaves of the strong-stable foliation $\mathcal{F}_1$ through these points, so we have $W^s(Q_k^-)\to W^{ss}(M^-)$ as $n\to +\infty$. Obviously, by an additional
small perturbation we can break the intersection between $W^u(P)$ and $W^{ss}(M^-)$ and create the heteroclinic intersection of
$W^u(P)$ with $W^s(Q_{k_0}^-)$, where $Q_{k_0}^-$ is some point from $\{Q_k^-\}$ (see figure \ref{fig:3d_3} (b)). 
\par{}
The last step is to show the existence of a transverse intersection of $W^s(P)$ and $W^u(Q_{k_0}^-)$. Denote by $P^+_{k^*}$ the point closest to $\Pi\cap W^s(O)$ of those points $P^+_k$ that survive the change of $\mu$. We will prove that $W^u(Q_{k_0}^-)$ intersects the stable manifold of $P^+_{k^*}$ using the expansion of two-dimensional areas by the Poincaré map. The non-empty intersection $W^s(P)\cap W^u(Q_{k_0}^-)$ follows from the homoclinic relation between $P$ and $P^+_{k^*}$. This completes the proof of the theorem.
\subsection{Heterodimensional cycles in a strange attractor}\label{sec:2.2}
Let us now consider the case where the above bifurcation happens within the strange attractor proposed in \citep{ts98}. In this paper, we show that heterodimensional cycles obtained by Theorems \ref{thm:hetero_1} can belong to such attractor and coexist there with a Newhouse wild set. Here a Newhouse wild set is a compact, invariant, and transitive hyperbolic set whose stable manifold intersects non-transversely its unstable manifold in a $C^2$ persistent fashion (see \citep{newhouse1,pv94,gts93}). 
\par{}
In order to have an attractor, we need the existence of a certain absorbing domain. We fix a neighbourhood of $O$ where
formula \eqref{eq:1setting_1} is valid; by a linear scaling of the variables, we can make the size of the neighbourhood equal to $1$. 
Let $S=\{(x,y_1,y_2,z)\mid|x|\leqslant 1,\|(y_1,y_2)\|=1,\|z\|\leqslant 1\}$ be a cross-section to $W^s_{loc}(O)$. Suppose that all orbits starting at $S$ return to $S$. Then, 
the region $D$ filled with all orbits of the flow starting from $S$, plus the equilibrium $O$ and its two unstable separatrices $\Gamma^+$ and $\Gamma^-$, is forward invariant. 
\par{}
As we mentioned before, the non-degeneracy condition (C1) imposed on the homoclinic loops  $\Gamma^+$ and $\Gamma^-$ along with the condition $\rho<1/2$ 
implies the volume hyperbolicity of the system near the set $O\cup \Gamma^+\cup \Gamma^-$. We now assume that this property extends to the whole
of the forward invariant region $D$:\\[10pt]
\indent (C5) (Volume-hyperbolicity condition): The tangent bundle of $D$ admits a continuous dominated splitting at any point of $D$: $T_D=N^{ss}\oplus N^c$, 
where $N^{ss}$ is the strong-stable subspace (corresponding to coordinates $z$ near $O$) and $N^c$ is the center subspace (corresponding to coordinates $x$ and $y$
near $O$). The flow restricted to $N^{ss}$ is exponentially contracting, and volumes are expanding in $N^c$.\\[10pt]
\indent Note that this condition implies the
existence of an absolutely-continuous invariant foliation tangent to $N^{ss}$ at each point of $D$ (see \citep{ts98}).
\par{}
Before we introduce the attractor, let us recall some definitions. Let $X_tP$ be the time shift of point $P$ by the flow $X$ for the time $t$. Take $\varepsilon>0$ and $\tau>0$; an $(\varepsilon,\tau)$-orbit is a sequence of points $P_1,P_2,\ldots ,P_k$ such that the distance between $P_{i+1}$ and $X_tP_i$ is smaller than $\varepsilon$ for some $t>\tau$. A point $Q$ is said to be 
$(\varepsilon,\tau)$-accessible from $P$ if there exists an $(\varepsilon,\tau)$-orbit connecting $P$ and $Q$, and accessible from $P$ if for some fixed $\tau$ and all $\varepsilon>0$, 
the point $Q$ is $(\varepsilon,\tau)$-accessible from $P$. A set $B$ is said to be accessible from a point $P$ if it contains a point that is accessible from $P$. A closed invariant set $B$ is called chain-transitive if, for any points $P$ and $Q$ in $B$ and for any $\varepsilon>0$ and $\tau>0$, the set $B$ contains an $(\varepsilon,\tau)$-orbit connecting $P$ and $Q$. A compact invariant set $B$ is called completely stable if, for any neighborhood $U(B)$, there exist $\varepsilon>0$ and $\tau>0$ and a neighbourhood $V(B)\subseteq U(B)$ such that all $(\varepsilon,\tau)$-orbits starting in $V(B)$ do not leave $U(B)$.
\par{}
The attractor of our system in $D$ is defined as the set $\mathcal{A}$ of all points accessible from $O$. It is shown in \citep{ts98} that $\mathcal{A}$ is the unique chain-transitive and 
completely stable set in $D$, and it is accessible from any point in $D$. Thus $\mathcal{A}$ is the unique Ruelle-Hurley attractor of the system in $D$ (see \citep{ru81,hu82}).
The volume-hyperbolicity implies that the attractor $\mathcal{A}$ is chaotic in the sense that every orbit in it has a positive maximal Lyapunov exponent. A complete
description of the structure of $\mathcal{A}$ is impossible, as it may contains a Newhouse wild-hyperbolic set (see \citep{ts98}). It is also shown in \citep{ts98} that
the attractor $\mathcal{A}$ may contain coexisting saddles of different indices. Here we strengthen the last statement by showing
the following corollary from Theorem \ref{thm:hetero_1}.
\begin{thm}\label{thm:hetero_2}
If system $X$ satisfies conditions (C2) - (C5), then there exists a system arbitrarily close to $X$ in $C^r$ such that it satisfies the same symmetry condition, and its attractor $\mathcal{A}$ in $D$ contains a Newhouse wild set, a symmetric pair of homoclinic loops, and
a symmetric pair of heterodimensional cycles near these loops.
\end{thm}
\par{}
As we mentioned in Section \ref{sec:1.2}, such strange attractor also exists in system \eqref{eq:example}.
%
%
%
%
%
%
%
%
%
%
\section{Proof of Theorem \ref{thm:hetero_1}}\label{sec:prf1}
We prove Theorem \ref{thm:hetero_1} by finding heterodimensional cycles for the Poincaré map on a cross-section $\Pi$ near the saddle-focus $O$. The proof is divided into several parts. We first describe the Poincaré map $T$ on $\Pi$. After this, we consider a two-parameter family $X_{\mu,\rho}$ with $X=X_{0,\rho^*}$. We find a sequence $\{(\mu_j,\rho_j)\}$ of parameter values accumulating on $(0,\rho^*)$ such that system $X_{\mu_j,\rho_j}$ has a homoclinic loop $\Gamma^-$ and a connection from a single-round index-2 periodic orbit to this loop (see Lemma \ref{lem:pm}). More specifically, we show that there exists an index-1 fixed point $P$ on the cross-section $\Pi$ such that $W^u(P)\cap W^{ss}(M^-)\neq\emptyset$ (where $M^-$ is the first intersection point of $\Gamma^-$ and $\Pi$). 
\par{}
Then, we define parameter $\nu$ which controls the separation of $W^u(P)$ and $W^{ss}(M^-)$ localised near the intersection point given by $(\mu_j,\rho_j)$. This means that, for any fixed pair $(\mu_j,\rho_j)$, the two manifolds $W^u(P)$ and $W^{ss}(M^-)$ will cross each other with a non-zero velocity as $\nu$ varies. Note that, in the space of dynamical systems with R-symmetry, there exists a codimension-1 surface $\mathcal{H}_1$ containing $X_{\mu_j,\rho_j}$ such that all systems on this surface have a double-round homoclinic loop $\Gamma^-$; inside this surface, there is a codimension-2 surface $\mathcal{H}_2$ corresponding to systems having the intersection $W^u(P)\cap W^{ss}(M^-)$. We now embed system $X_{\mu_j,\rho_j}$ into a two-parameter family $X_{\rho,\nu}$ such that this family lies in the surface $\mathcal{H}_1$ and systems in $\mathcal{H}_2$ correspond to parameter value $\nu=0$. We remark here that $\mu$ is now a function of $\rho$: when $\rho$ is changed, we need to change $\mu$ accordingly to keep the double-round homoclinic loop $\Gamma^-$. 
\par{}
Next, we obtain a sequence $\{(\rho_j^n,\nu_j^n)\}_n$ of parameter values converging to $(\rho_j,0)$ such that the unstable manifold $W^u(P)$ intersects the stable manifold of an index-2 periodic point $Q$ close to $M^-$ (see Lemma \ref{lem:quasisymmetric}). In the end, we show the existence of the non-empty intersection $W^s(P)\cap W^u(Q)$ for every pair $(\rho_j^n,\nu_j^n)$. Therefore, we obtain a heterodimensional cycle of the map $T$, which corresponds to one in the flow $X_{\rho_j^n,\nu_j^n}$. The second cycle exists by the symmetry.
\subsection{Construction of the Poincaré map $T$}\label{sec:prf1.1}
Recall that 
the local stable manifold $W^s_{loc}(O)$ is straightened and has the form $\{x=0\}$. We pick two points $M^+(0,y_1^+,0,z^+)$ and $M^-(0,y_1^-,0,z^-)$ near the equilibrium $O$ such that $M^+\in \Gamma^+\cap W^s_{loc}(O)$ and $M^-\in \Gamma^-\cap W^s_{loc}(O)$. We define $\Pi=\{(x,y_1,0,z)\mid|x|\leqslant\delta,|y_1-y^+|\leqslant\delta,\|z\|\leqslant\delta\}$ with an upper part $\Pi_1:=\Pi\cap\{ x> 0\}$ and a lower part $\Pi_2:=\Pi\cap\{ x< 0\}$. Denote by $\Pi_0$ the intersection of $\Pi$ with $W^s_{loc}(O)$, i.e. $\Pi\cap\{x=0\}$. Points on the cross-section have coordinates $(x,y_1,z)$ and we drop the subscript of $y_1$ for simplicity. Note that $\|z\|$ decreases much faster than $\|y\|$ along the homoclinic loops as $t\to+\infty$ so that we can assume $\|z^+\|<\delta$. Points $M^+$ and $M^-$ are the intersection points of $\Gamma^+$ and $\Gamma^-$ with $\Pi$, and have coordinates $(x^+,y^+,z^+)$ and $(x^-,y^-,z^-)$. By the symmetry, we have $x^+=-x^-$ and $y^+=y^-$. We will consider families of perturbed systems, so $x^\pm,y^\pm$ and $z^\pm$ are smooth functions of parameters. Since $\mu$ is the splitting parameter, we can assume $x^+=-x^-=\mu$
\par{}
In order to obtain the formula for the Poincaré map $T$, we use two additional cross-sections $\Pi_{glob_1}=\{(x=d,y_1,y_2,z)\}$ and $\Pi_{glob_2}=\{(x=-d,y_1,y_2,z)\}$, where $d>0$. The Poincaré map $T$ restricted to $\Pi_i$ $(i=1,2)$ is the composition of a local map $T_{loc_i}:\Pi_i\to \Pi_{glob_i}  , \,(x_0,y_0,z_0) \mapsto (y_1,y_2,z_1)$ and a global map $T_{glob_i}:\Pi_{glob_i}\to \Pi  , \,(y_1,y_2,z_1) \mapsto (\bar{x_0},\bar{y_0},\bar{z_0})$. We remark here that the ranges of $T_{glob_i}$ (and $T_i:=T_{glob_i}\circ T_{loc_i}$ mentioned later) are not completely contained in $\Pi$; however, when we iterate the Poincaré map later in the proofs, we will only consider the part which do not leave $\Pi$. Therefore, we write $T_{glob_i}:\Pi_{glob_i}\to \Pi$ (and $T_i:\Pi_i \to \Pi$ later) for simplicity.
\par{}
The map $T_{loc_i}$ is given by (see equation (13.4.13) of \citep{sstc2})
\begin{equation}\label{eq:setting_3}
\begin{array}{l}
 y_{1}=y_{0}\bigg(\dfrac{x_{0}}{(-1)^{i+1}d}\bigg)^{\rho }\cos \left( \omega \ln \left( \dfrac{(-1)^{i+1}d}{x_{0}}
\right) \right)+o\left( |x_{0}|^{\rho } \right) ,\\[20pt]
 y_{2}=y_{0}\bigg(\dfrac{x_{0}}{(-1)^{i+1}d}\bigg)^{\rho }\sin \left( \omega \ln \left( \dfrac{(-1)^{i+1}d}{x_{0}}
\right) \right)+o(|x_{0}|^{\rho }), \\[20pt]
z_1=\left(
\begin{array}{c}
o(|x_{0}|^{\rho })\\
\dots\\
o(|x_{0}|^{\rho })
\end{array} \right)
\end{array}
\end{equation}
\noindent Here we denote the small terms in each equation in \eqref{eq:setting_3} by $g_k$ $(k=1\dots n-1)$ , and we have
\begin{equation}\label{eq:smallterms}
\dfrac{\partial^{i+j} g_k}{\partial^i x \partial^j( y , z,\zeta,\rho)}=o(|x|^{\rho-i}) \quad i+j\leqslant(r-1).
\end{equation}
The global maps $T_{glob_i}$ are diffeomorphisms and can be written in Taylor expansions. We have
\begin{equation}\label{eq:setting_4}
\begin{array}{rcl}
 \bar{x}_{0}&=&(-1)^{i+1}(\mu +a_{11}y_{1}+a_{12}y_{2}+a_{13}z_1)+o\left(
|y_{1},y_{2},z_{1}| \right) \,, \\[10pt]
\bar{y}_{0}&=&y^+ +a_{21}x_{1}+a_{22}x_{2}+a_{23}z_1+o\left( |y_{1},y_{2},z_1| \right)\,, \\[10pt]
 \bar{z}_{0}&=&\mathcal{S}^{i+1}\left(z^+ +

 \begin{pmatrix}
 a_{31}x_{1}+a_{32}x_{2}+a_{33}z_1 \\
 \cdots \\
 a_{n-1,1}x_{1}+a_{n-1,2}x_{2}+a_{n-1,3}z_1
 \end{pmatrix} \right)
 
+o\left( |y_{1},y_{2},z_1| \right)\,,
 \end{array}
\end{equation}
\noindent where $a_{j3}$ $(j=1\dots n-1)$ are $(n-3)$-dimensional vectors. Let $T_1:=T|_{\Pi_1}=T_{glob_1}\circ T_{loc_1}:\Pi_1 \to \Pi$ and $T_2:=T|_{\Pi_2}=T_{glob_2}\circ T_{loc_2}:\Pi_2 \to \Pi$. Note that we have
\begin{equation}\label{eq:setting_6}
\lim_{M \to \Pi_0^+}T_1(M)=M^+ \quad  \mbox{and} \quad \lim_{M \to \Pi_0^-}T_2(M)=M^-.
\end{equation}
\noindent After the scalings $x={x_0}/{d},y={y_0}/{y^+}$ and $z=z_0$, and replacing ${\mu}/{d}$ by $\mu$, the maps $T_1$ and $T_2$ take the form
\begin{equation}\label{eq:1setting_4}
T_1:\quad\left\{\begin{array}{l}

                \bar{x}=\mu+Ayx^\rho\cos{(\omega \ln {\dfrac{1}{x}}+\theta)}+o(x^\rho), \\[15pt]
                \bar{y}=1+A_1yx^\rho\cos{(\omega \ln {\dfrac{1}{x}}+\theta_1)}+o(x^\rho), \\[15pt]
                \bar{z}=z^+ +
                \begin{pmatrix}
                 A_2yx^\rho\cos{(\omega \ln {\dfrac{1}{x}}+\theta_2)}+o(x^\rho) \\
                 \dots \\
                  A_{n-2}yx^\rho\cos(\omega \ln {\dfrac{1}{x}}+\theta_{n-2})+o(x^\rho) 
                \end{pmatrix},                
                \end{array}\right.
\end{equation}
\noindent and
\begin{equation}\label{eq:1setting_5}
T_2:\quad\left\{\begin{array}{l}
                \bar{x}=-\mu-Ay|x|^\rho\cos{(\omega \ln {\dfrac{1}{|x|}}+\theta)}+o(|x|^\rho), \\[15pt]
                \bar{y}=1+A_1y|x|^\rho\cos{(\omega \ln {\dfrac{1}{|x|}}+\theta_1)}+o(|x|^\rho), \\[15pt]
                \bar{z}=\mathcal{S} z^+ + \mathcal{S}
                \begin{pmatrix}
                 A_2y|x|^\rho\cos{(\omega \ln {\dfrac{1}{|x|}}+\theta_2)}+o(|x|^\rho) \\
                 \dots \\
                  A_{n-2}y|x|^\rho\cos(\omega \ln {\dfrac{1}{|x|}}+\theta_{n-2})+o(|x|^\rho) 
                \end{pmatrix} ,    
                 \end{array}\right.
\end{equation}
\noindent where $z \in \mathbb{R}^{n-3}$,   $A=y^+\sqrt{a_{11}^{2}+a_{12}^{2}}$, $
A_{1}=\sqrt{a_{21}^{2}+a_{22}^{2}}$, $ A_{m}=y^+\sqrt{a_{m+1,1}^{2}+a_{m+1,2}^{2}}$ $(m=2,\dots,n-2)$, ${\tan}\theta$ $
=-{a_{12}}/a_{11}$, $\tan\theta_{1}=-{a_{22}}/a_{21}$, $\tan\theta_{m}=-{a_{m+1,2}}/{a_{m+1,1}}$. Here $a_{ij}$ are coefficients from the global map given by formula \eqref{eq:setting_4}. The small terms $o(|x|^\rho)$ (for both $x>0$ and $x<0$) are functions of $x,y,z,\mu,\rho$ satisfying (\ref{eq:smallterms}). 
\par{}
Recall that we denote by $M^+$ and $M^-$ the first intersection points of $\Pi$ with $\Gamma^+$ and $\Gamma^-$. Their coordinates are now $(\mu,1,z^+)$ and $(-\mu,1,\mathcal{S}z^+)$. The maps $T_i$ $(i=1,2)$ and $T$ can be extended to $\Pi_i\cup\Pi_0$ and $\Pi$ respectively by letting
\begin{equation*}
T_1(0,y,z)=(\mu,1,z^+)\quad\mbox{and}\quad T_2(0,y,z)=(-\mu,1,\mathcal{S}z^+).
\end{equation*}
\par{}
From now on, we will work with the maps $T_1$ and $T_2$. The non-degeneracy condition mentioned in Section \ref{sec:1.1} is equivalent to
\begin{equation}\label{eq:pmn0}
AA_1\sin(\theta_1-\theta)\neq 0.
\end{equation}
\noindent Indeed, in the coordinate system satisfying (\ref{eq:1setting_1}) and (\ref{eq:1setting_2}), the transversality stated in the non-degeneracy condition is equivalent to the transversality of $T_{glob_1}(\Pi_{glob_1}$ $\cap W^{uE}_{loc}(O))$ and $T_{glob_2}(\Pi_{glob_2}\cap W^{uE}_{loc}(O))$ to the leaves $\{x=0,y=y^+\}$ through $M^+$ and $\{x=0,y=y^-\}$ through $M^-$, respectively, where the extended unstable manifold $W^{uE}_{loc}(O)$ is an invariant manifold tangent to the $\{z=0\}$ (see \citep{sstc1}).  By formula (\ref{eq:setting_4}), this is
\begin{equation*}
\mbox{det}\dfrac{\partial(\bar{x}_0,\bar{y}_0)}{(x_1,x_2)}\neq 0
\end{equation*} 
\noindent for both maps $T_{glob_1}$ and $T_{glob_2}$, which is equivalent to 
\begin{equation*}
\begin{array}{rcl}
\begin{vmatrix}
a_{11} & a_{12} \\
a_{21} & a_{22} \\
\end{vmatrix}=AA_1\sin(\theta_1-\theta)\neq 0
\end{array}.
\end{equation*}
\noindent As we mentioned before, if this condition is not satisfied, then we can make an arbitrary small perturbation to achieve it, and it will hold for all $C^r$-close systems.
\subsection{Coexistence of the homoclinic loop $\Gamma^-$ and an index-1 point $P$ with $W^u(P)\cap W^{ss}(M^-)\neq\emptyset$}\label{sec:prf1.2}
Let us now consider a two-parameter family $X_{\mu,\rho}$, where $X_{0,\rho^*}=X$. Shilnikov theorem implies that, for any system $X_{0,\rho}$ with $\rho<1$, there exists a countable set $\{P^+_k\}\subset\Pi_1$ of index-1 fixed points of $T_1$ accumulating on $M^+$. The proof of this theorem will be included in the proof of Lemma \ref{lem:WsP} in Section \ref{sec:prf1.5} where we show the non-empty transverse intersection $W^s(P)\cap W^u(Q)$. We now pick an arbitrary point $P$ from this set. In what follows, we consider sufficiently small perturbations such that $P$ remains an index-1 fixed point of $T_1$. The bound for the size of such perturbation will be given in Lemma \ref{lem:WsP}. We have the following result: 
\begin{lem}\label{lem:pm}
For any given value $\rho^* \in (0,{1}/{2})$, there exists a sequence $\{(\mu_j,\rho_j)\}$ accumulating on $(0,\rho^*)$ such that the corresponding system $X_{\mu_j,\rho_j}$ has a double-round homoclinic loop $\Gamma^-$ whose first intersection point $M^-$ with $\Pi$ has a strong-stable manifold $W^{ss}(M^-)$ that intersects the unstable manifold $W^u(P)$ quasi-transversely.
\end{lem}
Here {\it quasi-transversality} means that, for two manifolds $U$ and $V$, we have $\mathcal{T}_x U\cap\mathcal{T}_x V=\{0\}$ for the intersection point $x$ of $U\cap V$, where $\mathcal{T}_xU$ and $\mathcal{T}_xV$ are tangent spaces. The intersection $W^u(P)\cap W^{ss}(M^-)$ is quasi-transverse if it exists. Indeed, the strong stable manifold $W^{ss}(M^-)$ is a leaf of the foliation $\mathcal{F}_1$ tangent to strong-stable directions (i.e. close to $z$-directions), and $W^u_{loc}(P)$ is tangent to the center-unstable direction (i.e. close to $(x,y)$-directions). Therefore, by letting $\{M\}=W^u(P)\cap W^s(Q)$, we have $\mathcal{T}_M W^s_{loc}(Q) \cap \mathcal{T}_M W^u_{loc}(P)=\{0\}$, which gives the quasi-transversality.
\par{}
\noindent{\it Proof of Lemma \ref{lem:pm}}. We first change $\mu$ to make $\Gamma^-$ a double-round homoclinic loop. This can be done by solving the equation $T_2(M^-)=(0,y,z)$, where $M^-(-\mu,1,\mathcal{S}z^+)$ is the first intersection point of $\Gamma^-$ with $\Pi$, and $y,z$ can be arbitrary. By formulas (\ref{eq:1setting_4}) for $T_1$ and (\ref{eq:1setting_5}) for $T_2$, we have
\begin{equation*}\label{eq:pm1}
0=\mu+A|\mu|^{\rho}\cos(\omega\ln\dfrac{1}{|\mu|}+\theta)+o(|\mu|^{\rho})  \quad\mbox{if}\quad \mu<0,
\end{equation*}
\noindent and
\begin{equation*}\label{eq:pm2}
0=-\mu-A\mu^{\rho}\cos(\omega\ln\dfrac{1}{\mu}+\theta)+o(\mu^{\rho}) \quad\mbox{if}\quad \mu>0.
\end{equation*}
\noindent Denote
\begin{equation}\label{pmn1}
\omega\ln\dfrac{1}{|\mu|}=2\pi j_0+ \xi_0 -\theta \quad \xi_0\in[0,2\pi).
\end{equation}
\noindent The above two equations can be rewritten as
\begin{equation}\label{eq:pm3}
0=\exp\bigg(\dfrac{-2\pi j_0- \xi_0 +\theta}{\omega}\bigg)-A\exp\bigg(\dfrac{-2\pi j_0\rho- \xi_0\rho +\theta\rho}{\omega}\bigg)\cos\xi_0+o\bigg(\exp\bigg(\dfrac{-2\pi j_0\rho}{\omega}\bigg)\bigg),
\end{equation}
\noindent and
\begin{equation}\label{eq:pm4}
0=\exp\bigg(\dfrac{-2\pi j_0- \xi_0 +\theta}{\omega}\bigg)+A\exp\bigg(\dfrac{-2\pi j_0\rho- \xi_0\rho +\theta\rho}{\omega}\bigg)\cos\xi_0+o\bigg(\exp\bigg(\dfrac{-2\pi j_0\rho}{\omega}\bigg)\bigg),
\end{equation}
\noindent respectively. 
\par{} 
Since $\rho<1$ and $j_0$ is large, we have
\begin{equation*}
\exp\bigg(\dfrac{-2\pi j_0- \xi_0 +\theta}{\omega}\bigg)\ll\exp\bigg(\dfrac{-2\pi j_0\rho- \xi_0\rho +\theta\rho}{\omega}\bigg)
\end{equation*} 
We divide both sides of equations \eqref{eq:pm3} and \eqref{eq:pm4} by $\exp\big(\dfrac{-2\pi j_0\rho}{\omega}\big)$ and take the limit $j_0\to+\infty$. We seek for the solutions to the limit systems, which, by implicit function theorem, will give us solutions to the original systems. Either of equations \eqref{eq:pm3} and \eqref{eq:pm4} leads to
\begin{equation*}
\cos\xi_0=o(1)_{j\to +\infty}.
\end{equation*}
This along with equation \eqref{pmn1} gives a sequence $\{\mu_{j_0}\}$ of solutions of the form
\begin{equation}\label{pmn2}
\mu_{j_0}=\pm\exp\Bigg(\dfrac{-2\pi j_0-\dfrac{\pi}{2}-m\pi+\theta}{\omega}\Bigg)+o(1)_{j_0\to+\infty}
\quad m=0,1.
\end{equation}
\noindent Obviously, $\mu_{j_0}\to 0$ as $j_0\to+\infty$. Such values of $\mu$ give us a double-round homoclinic loop $\Gamma^-$ (and another one $\Gamma^+$ by symmetry).
\par{}
Let us now find the intersection of $W^u(P)$ with $W^{ss}(M^-)$. Denote the coordinates of $P$ by $(x_p,y_p,z_p)$. By taking a vertical line joining $P$ and a point on $\Pi_0$ and iterating it, one can check that the local unstable manifold $W^u_{loc}(P)$ of $P$ is spiral-like and winds onto $M^+$, which is given by
\begin{equation}\label{eq:pm5}
\begin{array}{rcl}
x&=&\mu+Ay_pt^{\rho}\cos(\omega\ln\dfrac{1}{t}+\theta)+o(t^\rho) ,\\[10pt]
y&=&1+A_1y_pt^{\rho}\cos(\omega\ln\dfrac{1}{t}+\theta_1)+o(t^\rho) , \\[10pt]
z&=&z^+ +
                \begin{pmatrix}
                 A_2y_p t^\rho\cos{(\omega \ln {\dfrac{1}{t}}+\theta_2)}+o(t^\rho) \\
                 \dots \\
                  A_{n-2}y_p t^\rho\cos(\omega \ln {\dfrac{1}{t}}+\theta_{n-2})+o(t^\rho) 
                \end{pmatrix} 
\end{array}
\end{equation}
\noindent where $t\in(0,x_p)$. 
\par{}
We need a formula for the local stable manifold $W^{s}_{loc}(Q_1)$, which is a leaf of the strong-stable foliation $\mathcal{F}_1$. The leaves of $\mathcal{F}_1$ are given by the following lemma.
\begin{lem}\label{lem:WsQ}
Let $M(x_0,y_0,z_0)$ be a point on $\Pi$ with $y_0$ sufficiently small. The local strong-stable manifold $W^{ss}_{loc}(M)$ (i.e. the leaf of $\mathcal{F}_1$ through $M$) is the graph of the function
\begin{equation*}
(x,y)=(x_0+(z-z_0)a_1,y_0+(z-z_0)a_2),
\end{equation*}
\noindent where $a_1=o(|x_0|)$ and $a_2=o(1)_{x_0\to0}$ are $(n-3)$-dimensional vectors whose components are certain functions of $z,y_0,x_0,z_0$ and the parameters (denoted by a vector $\varepsilon$) satisfying
\begin{equation*}
\dfrac{\partial^{i+l} a_1}{\partial^i x_0 \partial^l (z,y_0, z_0,\varepsilon)}=o(|x_0|^{1-i})\quad i+l\leqslant(r-1),
\end{equation*}
\noindent and
\begin{equation*}
\dfrac{\partial^{l} a_2}{\partial^l (z,y_0, z_0,\varepsilon)}=o(1)_{x_0\to 0} \quad l\leqslant(r-1) . 
\end{equation*}
\end{lem}
The proof of this lemma is postponed until Section \ref{sec:prf1.6}. The strong-stable manifold $W^{ss}(M^-)$ is now given by
\begin{equation}\label{eq:pm6}
\begin{array}{rcl}
x&=&-\mu+(z+z^+)o(|\mu|)  ,\\[10pt]
y&=&1+(z+z^+)o(1)  .
\end{array}
\end{equation}
\noindent The intersection point of $W^{ss}(M^-)$ with $W^u(P)$ is found by simultaneously solving equations (\ref{eq:pm5}) and (\ref{eq:pm6}). By noting $y_p=1+O(x_p^\rho)$ from formula \eqref{eq:1setting_4} and $z=z^++O(t^\rho)$ from \eqref{eq:pm5}, finding the intersection $W^{ss}(M^-)\cap W^u(P)$ is equivalent to solving the equations
\begin{equation}\label{eq:pm7}
\begin{array}{rcl}
2\mu&=&-At^\rho\cos(\omega\ln\dfrac{1}{t}+\theta)+o(t^\rho)+o(|\mu|),\\[10pt]
o(1)&=&A_1t^{\rho}\cos(\omega\ln\dfrac{1}{t}+\theta_1)+o(t^\rho) .
\end{array}
\end{equation}
\noindent This can now be seen as finding an intersection of a spiral given by the RHS of system \eqref{eq:pm7} of equations with a point $(u(\mu),v(\mu)):=(2\mu,o(1))$. Note that the $\mu$ value is given by equation \eqref{pmn2}, so here we will solve \eqref{eq:pm7} for $t$ and $\rho$. 
\par{}
We first find $t$. Let $\tan(\theta_1-\theta)=-{b}/{a}$ and rewrite equation (\ref{eq:pm7}) as
\begin{equation*}\label{eq:pm13}
\begin{array}{rcl}
u&=&-At^\rho\cos(\omega\ln\dfrac{1}{t}+\theta)+\dots,\\[10pt]
v&=&\dfrac{A_1}{\sqrt{a^2+b^2}}t^{\rho}\Big(b\cos(\omega\ln\dfrac{1}{t}+\theta)+a\sin(\omega\ln\dfrac{1}{t}+\theta)\Big)+\dots ,
\end{array}
\end{equation*} 
\noindent where we denote the small terms that tend to zero as $\mu,t$ tend to zero by dots throughout the proof. The above equations yield
\begin{equation*}
\dfrac{v}{u}=\dfrac{A_1}{-A\sqrt{a^2+b^2}}(b+\tan(\omega\ln(\dfrac{1}{t})+\theta))+\dots,
\end{equation*}
\noindent i.e.
\begin{equation}\label{eq:pm14}
\omega\ln t=\theta+\arctan\Big(\dfrac{Av\sqrt{a^2+b^2}}{A_1au}+\dfrac{b}{a}\Big)+2\pi k+\dots.
\end{equation}
\noindent Note that $u=2\mu$ does not vanish (see equation \eqref{pmn2}), so no matter how $\mu$ and $\rho$ change, equation (\ref{eq:pm14}) has a solution $t$ which depends continuously on all parameters for every fixed $k$.
\par{}
We proceed to find the values for $\rho$. By plugging the coefficients $A=y^+\sqrt{a_{11}^{2}+a_{12}^{2}}$, $
A_{1}=\sqrt{a_{21}^{2}+a_{22}^{2}}$, $ {\tan}\theta=-a_{12}/a_{11}$ and $\tan\theta_{1}=-a_{12}/a_{11}$ (see discussion under equation \eqref{eq:1setting_5}) into \eqref{eq:pm7}, we have
\begin{equation}\label{eq:pm8}
\begin{array}{rcl}
u&=&-y^+a_{11}t^\rho\cos(\omega\ln\dfrac{1}{t})-x^+a_{12}t^\rho\sin(\omega\ln\dfrac{1}{t})+\dots,\\[10pt]
v&=&a_{21}t^{\rho}\cos(\omega\ln\dfrac{1}{t})+a_{22}t^{\rho}\sin(\omega\ln\dfrac{1}{t})+\dots ,
\end{array}
\end{equation}
\noindent where dots denote small terms that tend to zero as $\mu,t$ tend to zero. By the non-degeneracy condition \eqref{eq:pmn0}, the matrix
\begin{equation*}\label{eq:pm9}
\begin{pmatrix}
-y^+a_{11} & -y^+a_{12} \\
a_{21} & a_{22}
\end{pmatrix}
\end{equation*}
\noindent is invertible. We denote the inverse matrix as 
\begin{equation*}\label{eq:pm10}
\begin{pmatrix}
b_{11} & b_{12} \\
b_{21} & b_{22}
\end{pmatrix}.
\end{equation*}
\noindent Then, we can rewrite equations in (\ref{eq:pm8}) as
\begin{equation*}\label{eq:pm11}
\begin{array}{rcl}
b_{11}u+b_{12}v&=&t^\rho\cos(\omega\ln\dfrac{1}{t})+\dots,\\[10pt]
b_{21}u+b_{22}v&=&t^{\rho}\sin(\omega\ln\dfrac{1}{t})+\dots ,
\end{array}
\end{equation*}
\noindent which, by squaring, summing up the above equations and taking logarithm, gives
\begin{equation}\label{eq:pm12}
\rho\ln t=\dfrac{1}{2}\ln(b_{11}^2u^2+b_{12}^2v^2+b_{21}^2u^2+b_{22}^2v^2+2(b_{11}b_{12}+b_{21}b_{22})uv)+\dots.
\end{equation} 
\par{}
Let us divide equation \eqref{eq:pm12} by \eqref{eq:pm14} and consider the ratio 
\begin{equation}\label{eq:pm15}
\dfrac{\rho}{\omega}=\dfrac{\ln(b_{11}^2u^2+b_{12}^2v^2+b_{21}^2u^2+b_{22}^2v^2+2(b_{11}b_{12}+b_{21}b_{22})uv)}{2\theta+2\arctan\Big(\dfrac{Av\sqrt{a^2+b^2}}{A_1au}+\dfrac{b}{a}\Big)+4\pi k}+\dots.
\end{equation} 
\noindent Note that $v$ may change as $\rho$ and $\omega$ change, and the above equation is not an explicit function for $\rho$. We will show that we can find values of $\rho$ from this equation anyway, and moreover they form a dense set when $j_0$ and $k$ tend to infinity. 
\par{}
The numerator of equation \eqref{eq:pm15} satisfies
\begin{equation}\label{eq:pm16}
C\mu^2<b_{11}^2u^2+b_{12}^2v^2+b_{21}^2u^2+b_{22}^2v^2+2(b_{11}b_{12}+b_{21}b_{22})uv<o(1)_{\mu\to 0},
\end{equation}
\noindent where $C$ is a constant independent of $\rho,\omega,\mu$ and $t$. Note that the coefficients $b_{ij}$ depend on all parameters. However, the range of parameter change is small, so the coefficients just vary slightly. This means that the constant $C$ can be chosen the same for all parameters under consideration, and it remains bounded away from zero and infinity. We denote the right hand side of equation (\ref{eq:pm15}) by $H$. Inequality \eqref{eq:pm16}, equation \eqref{pmn2} along with the fact that the value 
\begin{equation*}
\arctan\Big(\dfrac{Av\sqrt{a^2+b^2}}{A_1au}+\dfrac{b}{a}\Big)
\end{equation*}
\noindent is bounded imply that there exist two functions $K_1(j_0)$ and $K_2(j_0)$ such that 
\begin{equation}\label{eq:pm17}
\dfrac{K_1(j_0)}{k}<H<\dfrac{K_2 (j_0)}{k}.
\end{equation}
Here $K_1(j_0)$ and $K_2(j_0)$ do not depend on $t$ and parameters, and we have $K_1(j_0),K_2(j_0)\to +\infty$ as $j_0\to +\infty$. We now consider the function 
\begin{equation}\label{eq:pm18}
G(\rho)=\dfrac{\rho}{\omega}-H.
\end{equation}
By continuity, whatever $j_0$ and $k$ we choose, we can find a value of $\rho$ such that $G(\rho)=0$ by changing $\rho$ from ${K_1(j_0)}/{k}$ to ${K_2(j_0)}/{k}$. Note that, for any given number, we can choose a sequence of $(j_0,k)$ such that the corresponding intervals $({K_1(j_0)}/{k},{K_2(j_{0})}/{k})$ shrink to this number as $(j_0,k) \to (+\infty,+\infty)$. Hence, we obtain a dense set of $\rho$ values. Lemma \ref{lem:pm} is proven.\qed
\subsection{Countable sets of index-2 periodic points}\label{sec:prf1.3}
We first prove a lemma on the condition for a periodic point of $T$ to have index 2. Then we will show that, for some parameter values, there are infinitely many index-2 periodic points near the intersection points of $\Gamma^+$ and $\Gamma^-$ with $\Pi$.
\par{}
We introduce a transformation for the $x$-coordinates of points on the cross-section $\Pi$: 
\begin{equation}\label{eq:mr1.11}
\omega\ln\dfrac{1}{x}=2\pi j+ \xi -\theta  ,\quad\xi\in[0,2\pi),
\end{equation}
\noindent by which we divide the cross-section into different regions, and $\xi$ is a new coordinate in each region. Let $Q$ be a periodic point of $T$ of period $k$ and have the orbit $\{Q=Q_1(x_1,y_1,z_1),$ $Q_2(x_2,y_2,z_2),\dots,Q_k(x_k,y_k,z_k)\}$. The $x$-coordinates of this orbit are represented as
\begin{equation}\label{eq:index2_1}
\omega\ln\dfrac{1}{|x_i|}=2\pi j_i+ \xi_i -\theta \quad i=1,2,\dots,k.
\end{equation} 
\begin{lem}\label{lem:index_2}
The period-k point $Q$ is of index 2, if and only if
\begin{equation*}
\cos(\xi_1-\varphi)\cos(\xi_2-\varphi)\dots\cos(\xi_k-\varphi)=c\psi(\xi,j,y,z) \quad -1<c<1 ,
\end{equation*}
\noindent where $\varphi=\arctan({\omega}/{\rho})$, $\xi=(\xi_1,\xi_2,\dots,\xi_k)$, $j=(j_1,j_2,\dots,j_k)$, $y=(y_1,y_2,\dots,y_k)$, $z=(z_1,z_2,\dots,z_k)$ and $\psi=o(1)_{j_1,j_2,\dots,j_k \to \infty}$ is a certain function depending continuously on $\xi,y,z$ and parameters $\mu,\rho$, such that 
\begin{equation*}
\dfrac{\partial^{i} \psi}{\partial^i ({\xi}, {y} , {z}, \mu,\rho)}=o(1)_{j_1,j_2,\dots,j_k \to \infty} \quad (i\leqslant(r-2)).
\end{equation*}
\end{lem}
\noindent{\it Proof.} We start by computing the trace of the matrix $\mbox{D}T^{(k)}:={\partial T^{(k)}}/{\partial( x , y ,z)}$, which is the product of matrices of the form $\mbox{D}T_1$ or
$\mbox{D}T_2$ depending on the orbit of $Q$. Note that $\mbox{D}T_1$ and
$\mbox{D}T_2$ are the same up to different coefficients in front of functions of the coordinates in each entry. Therefore, as can be seen from the computation below, we obtain the same result for any type of composition of $\mbox{D}T^{(k)}$. For certainty, we assume that $Q$ is periodic under $T_1$.
\par{}
From formula \eqref{eq:1setting_4} for $T_1$, one can check that the $y$-coordinate of $Q_i$ $(i=1,\dots,k)$ satisfies $y_i=1+O(x_1 ^\rho, x_2^\rho,\dots , x_k^\rho)$. Thus, $y_i$ can be sufficiently close to 1 if we choose $x_1, x_2,\dots , x_k$ sufficiently close to zero. By transformation \eqref{eq:mr1.11} on coordinate $x$ and formula \eqref{eq:1setting_4}, we have
\begin{equation}\label{eq:WsQ_1}
\begin{array}{c}
\left.\mbox{D}T_1\right|_{Q_i}=\\[30pt]
\left(
\begin{array}{lll}
Ax_i ^{\rho-1}(\rho\cos\xi_i+\omega\sin\xi_i)+o(x_i ^{\rho-1}) & Ax_i ^\rho\cos\xi_i+o(x_i ^\rho) &\boldsymbol{a} \\[10pt]

\begin{array}{l}
-A_1x_i ^{\rho -1}(\rho \cos(\xi_i+\theta_1-\theta)\\
+\omega\sin(\xi_i+\theta_1-\theta))+o(x_i ^{\rho -1})
\end{array}
&
\begin{array}{l}
A_1x_i ^\rho\cos(\xi_i+\theta_1-\theta)\\
+o(x_i ^\rho)
\end{array}
 & \boldsymbol{a}_1\\[20pt]

\begin{array}{l}
-A_2x_i ^{\rho -1}(\rho \cos(\xi_i+\theta_2-\theta)\\
+\omega\sin(\xi_i+\theta_2-\theta))+o(x_i ^{\rho -1})
\end{array}
&
\begin{array}{l}
A_2x_i ^\rho\cos(\xi_i+\theta_2-\theta)\\
+o(x_i ^\rho)
\end{array}
 &\boldsymbol{a}_2\\[10pt]

\dots&\dots&\dots \\[10pt]

\begin{array}{l}
-A_{n-2}x_i ^{\rho -1}(\rho \cos(\xi_i+\theta_{n-2}-\theta)\\
+\omega\sin(\xi_i+\theta_{n-2}-\theta))+o(x_i ^{\rho -1})
\end{array}
&
\begin{array}{l}
A_{n-2}x_i ^\rho\cos(\xi_i+\theta_{n-2}-\theta)\\
+o(x_i ^\rho)
\end{array}
 &\boldsymbol{a}_{n-2}
 \\

\end{array}
\right),\\
\end{array}
\end{equation}\\[5pt]
\noindent where $\boldsymbol{a}$ and $\boldsymbol{a}_j$ $(j=1\dots n-2)$ are $n-3$ dimensional vectors (rows) of the form $(o(x_i ^\rho),\dots,o(x_i ^\rho))$. We symbolically represent matrix \eqref{eq:WsQ_1} in the form
\begin{equation*}
\begin{pmatrix}
a_{11}^i & a_{12}^i & \boldsymbol{a}_{13}^i \\
\dots & \dots & \dots \\
{a}_{n-1,1}^i & {a}_{n-1,2}^i & \boldsymbol{a}_{n-1,3}^i 
\end{pmatrix}.
\end{equation*}
\par{}
Now we have
\begin{equation}\label{eq:index2_2}
\begin{array}{rcl}
\left.\mbox{D}T_1^{(k)}\right|_{Q_k}&=&
\begin{pmatrix}
a_{11}^1 & a_{12}^1 & \boldsymbol{a}_{13}^1 \\
\dots & \dots & \dots \\
{a}_{n-1,1}^1 & {a}_{n-1,2}^1 & \boldsymbol{a}_{n-1,3}^1 
\end{pmatrix}
\begin{pmatrix}
a_{11}^2 & a_{12}^2 & \boldsymbol{a}_{13}^2 \\
\dots & \dots & \dots \\
{a}_{n-1,1}^2 & {a}_{n-1,2}^2 & \boldsymbol{a}_{n-1,3}^2 
\end{pmatrix}
\dots
\begin{pmatrix}
a_{11}^k & a_{12}^k & \boldsymbol{a}_{13}^k \\
\dots & \dots & \dots \\
{a}_{n-1,1}^k & {a}_{n-1,2}^k & \boldsymbol{a}_{n-1,3}^k 
\end{pmatrix}
\\[30pt]
&=&
\begin{pmatrix}
a_{11}^1 a_{11}^2\dots a_{11}^k+\dots & \dots & \dots\\
\dots & \dots & \dots \\
\dots & \dots & \dots 
\end{pmatrix}.
\end{array}
\end{equation}
\noindent Here each of the terms denoted by dots contains at least one of $x_i^\rho$ ($i=1\dots k$) as a factor, and the term $a_{11}^1 a_{11}^2\dots a_{11}^k$ is the only one which does not contain $x_i^\rho$ and therefore the dominant one. This gives us the formula for the trace:
\begin{equation}
\begin{split}\label{eq:index2_3}
\begin{array}{rcl}
\mathrm{tr}\left.\mbox{D}T_1^{(k)}\right|_{Q_i}&=&a_{11}^1 a_{11}^2\dots a_{11}^k+o(x_1 ^{\rho-1}+\dots+x_k^{\rho-1})\\
&=&A^k(\rho^2+\omega^2)^{\frac{k}{2}}x_1 ^{\rho-1}\dots x_k^{\rho-1}\cos{(\xi_1-\varphi)}\dots \cos{(\xi_k-\varphi)}  
+o(x_1 ^{\rho-1}+\dots+x_k^{\rho-1}),
\end{array}
\end{split}
\end{equation}
\noindent where $\varphi=\arctan{{\omega}/{\rho}}$.
\par{}
Before we proceed further, we show that the eigenvalues of $\left.\mbox{D}T_1^{(k)}\right|_{Q_i}$ corresponding to $z$ coordinates have the order of $o(x_1 ^{\rho}\dots x_k^\rho)$. We use Lemma \ref{lem:WsQ}, which implies that, for any periodic point, there exists an $(n-3)$-dimensional invariant subspace $\mathbb{E}_z$ to which its strong-stable manifold is tangent. Note that the eigenvalues of $\left.\mbox{D}T_1^{(k)}\right|_{Q_i}$ can be divided into two groups: one includes $\lambda_1$ and $\lambda_2$ corresponding to coordinates $x$ and $y$; the other includes $\lambda_3\dots\lambda_{n-1}$ corresponding to the restriction $S$ of $\left.\mbox{D}T_1^{(k)}\right|_{Q_i}$ to $\mathbb{E}_z$. Since $\mathbb{E}_z$ is in the cone defined in the proof of Lemma \ref{lem:WsQ}, we have $\|S\Delta z\|\leqslant o(x_1 ^{\rho}\dots x_k^\rho)\|\Delta z\|$ (see (\ref{eq:WsQ_7})), where $\Delta z\in\mathbb{E}_z$. It follows that the $(n-3)$ strong-stable multipliers $\lambda_3\dots\lambda_{n-1}$ satisfy $\lambda_i=o(x_1 ^{\rho}\dots x_k^\rho)$, where $i=3,4\dots n-1$. 
\par{}
We have the following expressions:
\begin{equation}\label{eq:index2_4}
\lambda_1+\lambda_2=\mathrm{tr}\left.\mbox{D}T_1^{(k)}\right|_{Q_1}-\sum_{i=3}^{n-1}\lambda_i
\quad \mbox{and} \quad
\lambda_1 \lambda_2=\sum_{ij}M_{ij}-\sum_{\begin{subarray}{c}
1\leqslant i<j\leqslant n-1 \\
(i,j)\neq(1,2)
\end{subarray}}\lambda_i\lambda_j  ,
\end{equation}
\noindent where $M_{ij}$ is the minor obtained by taking $i$-th and $j$-th rows, and $i$-th and $j$-th columns from $\left.\mbox{D}T_1^{(k)}\right|_{Q_1}$. One can check that $M_{12}$ gives the largest contribution to $\sum_{ij}M_{ij}$ and that it is given by 
\begin{equation}\label{eq:index2_3a}
M_{12}=\prod_{i=1}^k{\left.\mbox{D}T_1\right|_{Q_i}}=C_{12}^k x_1 ^{2\rho -1}\dots x_k^{2\rho -1}+o(x_1 ^{2\rho -1}+\dots+ x_k^{2\rho -1}),
\end{equation}
\noindent where $C_{12}=-\omega AA_1\sin(\theta_1-\theta)$. Now from equations (\ref{eq:index2_3}) and (\ref{eq:index2_3a}), we have
\begin{equation}\label{eq:index2_5}
\begin{array}{rcl}
\lambda_1+\lambda_2&=&A^k(\rho^2+\omega^2)^{\frac{k}{2}}x_1 ^{\rho-1}\dots x_k^{\rho-1}\cos{(\xi_1-\varphi)}\dots\cos{(\xi_k-\varphi)}\\[10pt]
&&+o(x_1 ^{\rho-1}+\dots+x_k^{\rho-1})+o(x_1 ^{\rho}\dots x_k^\rho)
\end{array}
\end{equation}
\noindent and 
\begin{equation*}
\begin{array}{rcl}
\lambda_1 \lambda_2&=&C_{12}^k x_1 ^{2\rho -1}\dots x_k^{2\rho -1}+o(x_1 ^{2\rho -1}+\dots+x_k^{2\rho -1})\\[10pt]
&&+\lambda_1 o(x_1 ^{\rho}\dots x_k^\rho)+\lambda_2 o(x_1 ^{\rho}\dots x_k^\rho)+ o(x_1 ^{2\rho}\dots x_k^{2\rho}).
\end{array}
\end{equation*}
\noindent By noting $\lambda_1,\lambda_2<\|\mbox{D}T^k\|=O(x_1^{2\rho-1}\dots x_k^{2\rho-1})$, we have
\begin{equation}\label{eq:index2_6}
\lambda_1 \lambda_2=C_{12}^k x_1 ^{2\rho -1}\dots x_k^{2\rho -1}+o(x_1 ^{2\rho -1}+\dots+x_k^{2\rho -1}).
\end{equation} 
\par
\begin{figure}[!h]
\begin{center}
\includegraphics[scale=0.35]{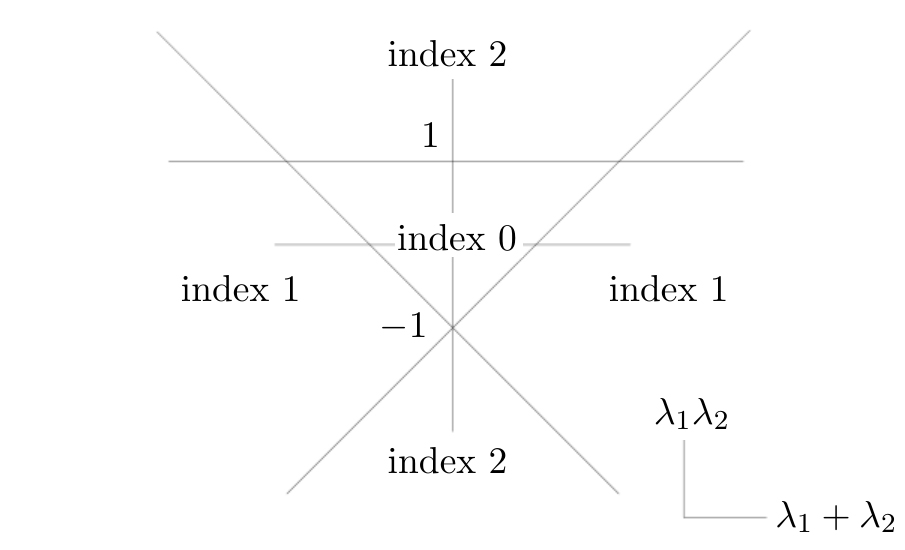}
\end{center}
\caption{The index is determined by the pair ($\lambda_1+\lambda_2,\lambda_1\lambda_2$).}
\label{pic:saddle region}
\end{figure}

{}
It can be checked that the equation for a periodic point to be of index 2 is
\begin{equation}\label{eq:index2_7}
\lambda_1+\lambda_2=c(\lambda_1\lambda_2+1), \quad c\in(-1,1).
\end{equation}
Indeed, the relation between the index of a point and the pair ($\lambda_1+\lambda_2,\lambda_1\lambda_2$) is illustrated in Figure \ref{pic:saddle region}, where we assume that the directions corresponding to $\lambda_3\cdots\lambda_{n-1}$ are the strong-stable ones.
\par{}
By noting relations \eqref{eq:index2_1} and dividing both sides of equation (\ref{eq:index2_7}) by $x_1 ^{\rho-1}x_2^{\rho-1}\dots$ $x_k^{\rho-1}$, we obtain
\begin{equation}\label{eq:index2_9}
\cos(\xi_1-\varphi)\cos(\xi_2-\varphi)\dots\cos(\xi_k-\varphi)=c\psi(\xi,y,z,j_1,j_2),
\end{equation}
\noindent where $\xi=(\xi_1,\xi_2,\dots,\xi_k)$, $\psi=o(1)_{j_1,j_2 \to \infty}$ is continuous in $\xi,y,z$ and all parameters such that
\begin{equation*}
\dfrac{\partial^{i} \psi}{\partial^i( {\xi} ,{y},{z},\zeta,\rho)}=o(1)_{j_1,j_2 \to \infty} \quad i\leqslant(r-2).
\end{equation*}
Since the computation of $\psi$ involves the first derivatives, we have $i\leqslant(r-2)$ instead of $i\leqslant(r-1)$.\qed
\par{}
Lemma \ref{lem:index_2} offers a criterion to seek for index-2 periodic orbits, by using which we have the following result.
\begin{lem}\label{lem:dense_i2}
For a dense set of $\rho$ values in $(0,\dfrac{1}{2})$, the system $X_{0,\rho}$ has two countable sets of index-2 periodic points of $T$ in $\Pi$ accumulating on $M^+$ and $M^-$, respectively.
\end{lem}
It can be seen from the proof below that finding each set of index-2 points only involves one homoclinic loop. Therefore, Lemma \ref{lem:dense_i2} extends Theorem 3 of \citep{os87} to systems with dimension higher than three, and we can reformulate it as follows.
\begin{cor}\label{cor:dense}
Let system $X_\rho$ be a $C^r$ flow in $\mathbb{R}^n$ ($r\geqslant3,n\geqslant4$) having a saddle-focus equilibrium $O$ with a homoclinic loop $\Gamma$ associated to it, and exactly one positive characteristic exponent with the ratio $\rho\in(0,{1}/{2})$. For a dense set of $\rho$ values in $(0,\dfrac{1}{2})$, system $X_{\rho}$ has a countable set of index-3 periodic orbits accumulating on $\Gamma$.
\end{cor}
\par{}
Note that the $\rho$ values obtained by Lemma \ref{lem:dense_i2} are not the same as those given by Lemma \ref{lem:pm}. We will show a way to obtain the quasi-transverse intersection $W^u(P)\cap W^s(Q)$ by using Lemmas \ref{lem:pm} and \ref{lem:dense_i2} together in the next subsection.
\par{}\noindent{\it Proof of Lemma \ref{lem:dense_i2}.} 
By symmetry it is sufficient to only consider $\Gamma^+$ and find a sequence of index-2 point accumulating on $M^+$.
\par{}
Let $\{Q_1(x_1,y_1,z_1), Q_2(x_2,y_2,z_2)\}\subset\Pi_1$ be an orbit of period 2 and index 2 under $T_1$. By \eqref{eq:mr1.11}, we have
\begin{equation}\label{eq:dense1}
\omega\ln\dfrac{1}{x_i}=2\pi j_i+ \xi_i -\theta  ,\quad\xi_i\in[0,2\pi)\quad i=1,2.
\end{equation}
\noindent By Lemma \ref{lem:index_2} and formula (\ref{eq:1setting_4}) for $T_1$, an orbit of $T_1$ with period 2 and index 2 is given by
{\allowdisplaybreaks
\begin{eqnarray*}\label{eq:dense2}
\qquad\qquad  &x_2=Ay_1x_1^\rho\cos{\xi_1}+o(x_1^\rho), \\[10pt]
               & y_2=1+A_1y_1x_1^\rho\cos{(\xi_1+\theta_1-\theta)}+o(x_1^\rho) ,\\[10pt]
                &z_2=z^+ +
                \begin{pmatrix}
                 A_2y_1x_1^\rho\cos{(\xi_1+\theta_2-\theta)}+o(x_1^\rho) \\
                 \dots \\
                  A_{n-2}y_1x_1^\rho\cos(\xi_1+\theta_{n-2}-\theta)+o(x_1^\rho) 
                \end{pmatrix}      ,\\[10pt]
               & x_1=Ay_2x_2^\rho\cos{\xi_2}+o(x_2^\rho), \\[10pt]
              &  y_1=1+A_1y_2x_2^\rho\cos{(\xi_2+\theta_1)}+o(x_2^\rho) ,\\[10pt]
              &  z_1=z^+ +
                \begin{pmatrix}
                 A_2y_2x_2^\rho\cos{(\xi_2+\theta_2-\theta)}+o(x_2^\rho) \\
                 \dots \\
                  A_{n-2}y_2x_2^\rho\cos(\xi_2+\theta_{n-2}-\theta)+o(x_2^\rho) 
                \end{pmatrix} ,\\[10pt]
               & \cos(\xi_1-\varphi)\cos(\xi_2-\varphi)=c\psi,
\end{eqnarray*}}

\noindent where $-1<c<1$ and $\psi$ is a certain function of $\xi,y,z$ depending continuously on parameters and $\psi\to 0$ as $j_1,j_2\to +\infty$. By expressing $y$ and $z$ as functions of $x$, we obtain a reduced system given by 
{\allowdisplaybreaks
\begin{eqnarray*}
&x_2=-A|x_1|^\rho\cos{\xi_1}+o(|x_1|^\rho)+O(|x_1|^\rho |x_2|^\rho), \label{eq:dense2.1}\\[10pt]            
               & x_1=-A|x_2|^\rho\cos{\xi_2}+o(|x_2|^\rho)+O(|x_1|^\rho |x_2|^\rho), \label{eq:dense2.2}\\[10pt]             
               & \cos(\xi_1-\varphi)\cos(\xi_2-\varphi)=c\psi. \label{eq:dense2.3}
\end{eqnarray*}}
It will be shown later that the following relation is satisfied:
\begin{equation}\label{eq:dense2.4}
|x_1|^\rho\sim |x_2|.
\end{equation}
\noindent Consequently, we replace the terms $o(|x_1|^\rho)+O(|x_1|^\rho |x_2|^\rho)$ and $o(|x_2|^\rho)+O(|x_1|^\rho |x_2|^\rho))$ in above equations by $o(|x_1|^\rho)$ and $o(|x_2|^\rho)$, respectively. The relation \eqref{eq:dense1} now brings these equations to the following form: 
\begin{equation}\label{eq:dense3}
\exp\bigg(\dfrac{-2\pi j_2-\xi_2+\theta}{\omega}\bigg)=A\exp\bigg(\dfrac{-2\pi\rho j_1-\rho\xi_1+\rho\theta}{\omega}\bigg)\cos\xi_1+o\bigg(\exp\bigg(\dfrac{-2\pi\rho j_1}{\omega}\bigg)\bigg) ,
\end{equation}\\[-0.5cm]
\begin{equation}\label{eq:dense4}
\exp\bigg(\dfrac{-2\pi j_1-\xi_1+\theta}{\omega}\bigg)=A\exp\bigg(\dfrac{-2\pi\rho j_2-\rho\xi_2+\rho\theta}{\omega}\bigg)\cos\xi_2+o\bigg(\exp\bigg(\dfrac{-2\pi\rho j_2}{\omega}\bigg)\bigg) ,
\end{equation}
\begin{equation}\label{eq:dense5}
\cos(\xi_1-\varphi)\cos(\xi_2-\varphi)=c\psi.
\end{equation}
\par{}
We solve this system with sufficiently large $j_1$ and $j_2$. One should note that the number of equations is larger than that of variables $(\xi_1,\xi_2)$, so whether this system of equations is solvable depends on the value of the parameter $\rho$. Throughout the rest of the proof, we denote by dots the small terms which are functions of $\xi_1,\xi_2,j_1,j_2$ and tend to zero as $j_1$ and $j_2$ tend to positive infinity. 
\par{}
Equation (\ref{eq:dense5}) implies that one of the terms $\cos(\xi_1-\varphi)$ and $\cos(\xi_2-\varphi)$ must be small. We assume that $\cos(\xi_1-\varphi)$ is small and $\cos(\xi_2-\varphi)$ is bounded away from zero, by which we have
\begin{equation}\label{eq:dense_5.5}
\cos(\xi_1-\varphi)=c\psi_1 ,
\end{equation}
\noindent where $\psi_1=o(1)_{j_1,j_2\to +\infty}$ is a function of $\xi_1,\xi_2,j_1,j_2$ depending continuously on $\xi_1,\xi_2$ and parameters. Then, we obtain
\begin{equation}\label{eq:dense5.1}
\xi_1=\arccos(c\psi_1)+k_1\pi+\varphi=\dfrac{\pi}{2}+k_1\pi+\varphi+\psi_2(\xi_1,\xi_2,j_1,j_2,c) ,
\end{equation}
\noindent where $k_1=0,1$ since $\xi_1\in[0,2\pi)$, $\varphi=\arctan({\rho}/{\omega})$, and the function $\psi_2=o(1)_{j_1,j_2\to +\infty}$ depends continuously on $\xi_1,\xi_2,c$ and parameters. Note that the value of $\psi_2$ changes slightly when $c$ varies in (-1,1).
\par{}
Another expression for $\cos\xi_1$ can by found from equation (\ref{eq:dense3}):
\begin{equation}\label{eq:dense_6}
\cos\xi_1=B^{-1}\exp\bigg(\dfrac{2\pi(\rho j_1-j_2)+\theta-\rho\theta+\rho\xi_1-\xi_2}{\omega}\bigg)+\dots.
\end{equation}
\noindent Recall that we assume that $\cos(\xi_1-\varphi)$ is small which means that $\cos\xi_1$ is bounded away from zero. Equation (\ref{eq:dense_6}) implies that $\cos(\xi_1)$ is positive. Therefore, we have $k_1=1$ in (\ref{eq:dense5.1}). 
\par{}
We proceed to find $\xi_2$. After taking logarithm on both sides of equation (\ref{eq:dense_6}) and sorting the terms, we get
\begin{equation}\label{eq:dense5.3}
\rho j_1-j_2=\omega\ln(B\cos\xi_1)-\theta+\rho\theta-\rho\xi_1+\xi_2+\dots,
\end{equation}
\noindent Since $\cos\xi_1$ is bounded away from zero, the above equation implies that $\rho j_1-j_2$ is bounded, which means $\rho j_2-j_1$ is large.
\par{}
We divide both sides of equation (\ref{eq:dense4}) by $\exp({-2\pi\rho j_2}/{\omega})$ and take the limit $j_1,j_2\to+\infty$. By noting that $\rho j_2-j_1$ is large, we have
\begin{equation}\label{eq:dense5.2}
\xi_2=\dfrac{\pi}{2}+k_2\pi+\dots ,
\end{equation}
\noindent where $k_2=0,1$ since $\xi_2\in[0,2\pi)$. The expression for $\xi_2$ implies that $\cos(\xi_2-\varphi)$ is bounded away from zero. This agrees with our assumption above \eqref{eq:dense_5.5}.
\par{}
The variables $\xi_i$ can be solved out by implicit function theorem from \eqref{eq:dense5.1} and \eqref{eq:dense5.2}, and they are functions of $j_1,j_2$ and $c$. By plugging the new expressions for $\xi_i$ into equation \eqref{eq:dense5.3}, we have
\begin{equation}\label{eq:dense5.4}
\rho j_1 - j_2=\Psi(\rho,j_1,j_2,c),
\end{equation}
\noindent where $c\in(-1,1)$ and $\Psi$ is a uniformly bounded function continuous in $\rho$ and $c$. For any fixed $\rho\in(1/2)$ and $c\in(0,1)$, each solution $(j_1,j_2)$ to equation \eqref{eq:dense5.4} gives an index-2 periodic orbit $\{Q_1,Q_2\}$ of the Poincaré map, and thus an index-3 periodic orbit of the system $X_{0,\rho}$. 
\noindent In what follows we show that there exists a dense set of $\rho$ values in $(0,1/2)$ such that, for each value $\rho$ in this set, there exists a sequence $\{j_1^n,j_2^n\}$ of solutions satisfying $j_2^n/j_1^n \to \rho$ as $n\to +\infty$. This will immediately imply Lemma \ref{lem:dense_i2}.
\par{}
By the boundedness of $\Psi(\rho,j_1,j_2,c)$, we can assume $|\Psi|<C$ for some constant $C$. Let $N_1>0$ be any large integer and $I_1=[a_1,a_2]$ be an arbitrary interval in $(0,1/2)$. We fix $c=c'\in(-1,1)$ and pick $j^1_1>N_1$. We have
\begin{equation}\label{eq:dense12}
a_1 j_1^1 - \Psi(\rho,j_1,j_2,c')<\lceil a_1 j_1^1 + C\rceil = m_1,
\end{equation}  
and
\begin{equation}\label{eq:dense13}
a_2 j_1^1 - \Psi(\rho,j_1,j_2,c')>\lfloor a_2 j_1^1 - C\rfloor = m_2.
\end{equation}
\noindent Note that $j_1^1$ can be chosen sufficiently large such that $m_2>m_1$. Now consider the function
\begin{equation*}
F(\rho,j_1,j_2):=\rho j_1 - \Psi(\rho,j_1,j_2,c') - m_1.
\end{equation*}
\noindent Equations \eqref{eq:dense12} and \eqref{eq:dense13} imply that 
\begin{equation}\label{eq:dense14}
F(a_1,j_1^1,m_1)<0 \quad\mbox{and}\quad F(a_2,j_1^1,m_1)>0.
\end{equation}
\noindent Therefore, by the continuity of $\Psi$, there exist a value $\rho_1\in I_1$ and a pair $(j_1^1,j_2^1=m_1)$ such that they satisfy
\begin{equation}\label{eq:dense15}
\rho_1 j_1^1 -j_2^1 = \Psi(\rho_1,j_1^1,j_2^1,c').
\end{equation}
\noindent
Note that $\Psi$ is also continuous in $c$, so one can find $c_1^1$ and $c_2^1$ with $c'\in(c_1^1,c_2^1)$ (or $(c_2^1,c_1^1)$) such that 
\begin{equation}\label{eq:dense16}
\Psi(\rho_1,j_1^1,j_2^1,c_1^1)<\rho_1 j_1^1 -j_2^1 <\Psi(\rho_1,j_1^1,j_2^1,c_2^1).
\end{equation}
The continuity of $\Psi$ in $\rho$ now implies that there exists a neighbourhood $J_1$ of $\rho_1$ such that inequality \eqref{eq:dense16} holds for all $\rho$ values taken from $J_1$.
\par{}
We choose $I_2\subset J_1$ such that $\rho_1\notin I_2$. We can find $\rho_2\in I_2$ and the corresponding pair $(c_1^2,c_2^2)$ with $c_1^2,c_2^2\in(0,1)$ and $(j_1^2,j_2^2)$ with $j_1^2>N_2>j_1^1$ such that an inequality of the same form of \eqref{eq:dense16} holds. By proceeding like this, we will find a sequence $\{I_n\}$ of nested intervals. Consequently, there exist a value $\rho_0\in \bigcap_{n=1}^{+\infty} I_n$, and two sequences $\{(c_1^n,c_2^n)\}$ and $\{(j_1^n,j_2^n)\}$ where $j_1^n,j_2^n\to +\infty$ and $j_2^n/j_1^n \to\rho_0$ as $n\to +\infty$ such that 
\begin{equation}\label{eq:dense17}
\Psi(\rho_0,j_1^n,j_2^n,c_1^n)<\rho_1 j_1^1 -j_2^1 <\Psi(\rho_0,j_1^n,j_2^n,c_2^n).
\end{equation}
This means that for each $n$ we can find a value $c^n\in(0,1)$ such that
\begin{equation}
\rho_0 j_1^n -j_2^n = \Psi(\rho_0,j_1^n,j_2^n,c^n),
\end{equation}
which implies the existence of an index-2 period-2 point of $T_2$.
\par{}
Since $I_1$ is chosen arbitrarily, such values $\rho_0$ are dense in $(0,{1}/{2})$.\qed
\subsection{Quasi-transverse intersection $W^u(P)\cap W^s(Q)$}\label{sec:prf1.4}
We fix a pair $(\mu_j,\rho_j)$ given by Lemma \ref{lem:pm} with a sufficiently large $j$, and consider perturbations of system $X_{\mu_j,\rho_j}$. Recall that parameter $\mu$ controls the splitting of the homoclinic loops $\Gamma^\pm$ around the points $\Gamma^\pm\cap\Pi$, and parameter $\rho$ is a function of the characteristic exponents at $O$. We now introduce another parameter $\nu$ which governs the separation of $W^{ss}(M^-)$ and $W^u(P)$ in system $X_{\mu_j,\rho_j}$ around the intersection given by Lemma \ref{lem:pm}. 
\par{}
To be more precise, we note that system $X_{\mu_j,\rho_j}$ lies in the codimension-1 surface $\mathcal{H}_1$ consisting of systems having the double-round homoclinic loop $\Gamma^-$. Inside $\mathcal{H}_1$, there exists a codimension-2 surface $\mathcal{H}_2$ such that the corresponding systems have the non-empty intersection $W^{ss}(M^-)\cap W^u(P)$. Note that the parameter $\mu$ is now a function of $\rho$ such that the system $X_{\rho,\nu}$ always has the double-round loop when we change $\rho$. These surfaces are smooth submanifolds of the space of sufficiently smooth $R-symmetric$????? systems (as we discussed under Theorem \ref{thm:hetero_1}, we make the original system $C^\infty$ from the very beginning). Indeed, if the system is $C^\infty$ then the invariant manifolds $W^s(O),W^u(O),W^u(P)$ and $W^{ss}(M^-)$ smoothly depend on the system. Obviously, this is true for $W^s(O),W^u(O)$ and $W^u(P)$ since the stable and unstable manifolds of hyperbolic periodic orbits and equilibrium states depend on parameters smoothly. Regarding $W^{ss}(M^-)$, we explain as follows. Inside the surface $\mathcal{H}_1$, the separatrix $\Gamma^-$ is always a homoclinic loop, and therefore, the strong-stable manifold $W^{ss}(M^-)$ lies in the stable manifold $W^s(O)$. By noting that $W^{ss}(M^-)$ is a leaf of the strong-stable foliation $\mathcal{F}_1$ and the part of $\mathcal{F}_1$ contained in $W^s(O)$ smoothly depends on parameters, we have that $W^{ss}(M^-)$ smoothly depends on the system in $\mathcal{H}_1$. As the corresponding stable and unstable manifolds depend on the system smoothly, the existence of their intersections correspond to the vanishing of certain smooth functionals. We can control the values of such functionals by adding to the system perturbations supported in a sufficiently small neighbourhood of the intersection points. This implies the smoothness of the surfaces $\mathcal{H}_1$ and $\mathcal{H}_2$.
\par{}
Let $\nu$ be the value of the functional that measures the splitting of $W^u(P)$ and $W^{ss}(M^-)$ for systems in $\mathcal{H}_1$, so the surface $\mathcal{H}_2$ is given by equation $\nu=0$. The explicit expression for $\nu$ is not obtained at this moment. It could be found if we had a better formula for the small terms $o(|x|^\rho)$ in (\ref{eq:1setting_4}) and (\ref{eq:1setting_5}). We embed the system $X^*:=X_{\mu_j,\rho_j}$ into a smooth two-parameter family $X_{\rho,\nu}$ of systems in $\mathcal{H}_1$. We assume that $X_{\rho,\nu}$ is transverse to $\mathcal{H}_2$ in $\mathcal{H}_1$. We also assume that the difference between $X_{\rho,\nu}$ and $X^*$ is localised in a small neighbourhood of $M^\pm,O$ and the preimage by $T$ of the intersection point of $W^u(P)\cap W^{ss}(M^-)$ given by Lemma \ref{lem:pm}. Therefore, changes of $\nu$ separate $W^u(P)$ and $W^{ss}(M^-)$ without changing the behaviour near the double-round homoclinic loop $\Gamma^-$, while changing $\rho$ is done without destroying the loop and the intersection of $W^u(P)$ with $W^{ss}(M^-)$. 
\par{}
Now we can apply Lemma \ref{lem:dense_i2} to a neighbourhood of the double-round homoclinic loop while keeping the intersection $W^{ss}(M^-)\cap W^u(P)$. This means that, by changing $\rho$ in the family $X_{\rho,\nu=0}$, one can find a value of $\rho$ arbitrarily close to $\rho_j$ such that the corresponding system will have a countable set $\{Q_k^-\}$ of periodic orbits of index 2 accumulating on $M^-$ while the intersection $W^u(P)\cap W^{ss}(M^-)$ is still intact. The points $Q_k^-$ have period 2 with respect to the double-round homoclinic loop $\Gamma^-$, so they have period 4 with respect to the original homoclinic loops. Since the stable manifolds $W^s(Q^-_k)$ and the strong-stable manifold $W^{ss}(M^-)$ are leaves of the foliation $\mathcal{F}_1$, we have that $W^s(Q^-_k)$ accumulates on $W^{ss}(M^-)$ as well. This sequence of index-2 points will persist when we change $\nu$ while keeping $\rho$ constant (as our family $X_{\rho,\nu}$ is such that the behaviour near the loop does not depend on $\nu$. Consequently, by an arbitrarily small change in $\nu$, we can destroy the intersection of $W^u(P)$ with $W^{ss}(M^-)$ and create an intersection of $W^u(P)$ with the stable manifold $W^{s}(Q)$ of a point $Q\in\{Q_k^-\}$. The quasi-transversality of this intersection $W^s(Q)\cap W^u(P)$ follows from the argument used above the proof of Lemma \ref{lem:pm}. Thus, we have the following result.

\begin{lem}\label{lem:quasisymmetric}
There exists a sequence $\{(\rho_j^n,\nu_j^n)\}_n$ where $(\rho_j^n,\nu_j^n)\to(\rho_j,0)$ as $n\to+\infty$ such that the separatrix $\Gamma^-$ forms a double-round homoclinic loop, and, on the cross-section $\Pi$ in system $X_{\rho_j^n,\nu_j^n}$, there exits an index-1 fixed point $P$ and an index-2 period-4 point $Q$ for which the intersection $W^u(P)\cap W^s(Q)$ is non-empty. Moreover, this intersection is quasi-transverse.
\end{lem}
\subsection{Transverse intersection $W^s(P)\cap W^u(Q)$}\label{sec:prf1.5}
As mentioned before, by Shilnikov theorem, there exist two countable sets $\{P_k^+\}\subset\Pi_1$ and $\{P_k^-\}\subset\Pi_2$ of index-1 fixed points of $T_1$ and $T_2$ accumulating on $M^+$ and $M^-$, respectively. The points $P_k^+$ and $P_k^-$ are obtained by solving the equations $T_1(x,y,z)=(x,y,z)$ and $T_2(x,y,z)=(x,y,z)$, respectively. In the previous sections, we picked an arbitrary point $P$ from the set $\{P_k^+\}$ which remains an index-1 fixed point under the small perturbation. Note that all results still hold if we pick $P$ from the set $\{P_k^-\}$. In what follows, we first give a detailed discussion on the points $P_k^+$ and their local stable manifolds; then we prove that the intersection $W^s(P)\cap W^u(Q)$ is non-empty, 
\subsubsection{The set $\{P_k^+\}$ of index-1 fixed points}\label{sec:prf1.5.1}
We start by giving a lemma on the local stable manifolds of the points $P^+_k$ in two cases, where we consider the $\mu=0$ case for the original set of index-1 fixed points given by Shilnikov theorem, as well as the $\mu\neq 0$ case for those points which survive the perturbation.
\begin{lem}\label{lem:WsP}
At $\mu=0$, there exists a set $\{P_k^+\}$ of index-1 fixed points of the map $T_1$ such that the $x$-coordinate of each point is given by
\begin{equation}\label{eq:WsP0}
x_k=C\exp\Big(\dfrac{-\pi k}{\omega}\Big)+o\Big(\exp\Big(\dfrac{-\pi k}{\omega}\Big)\Big),
\end{equation}
\noindent where $C=\exp(({2\theta-\pi})/{2\omega})$. The local stable manifolds $W^s_{loc}(P_k^+)$ are graphs of functions $g(y,z)$ defined for all $y$ and $z$ values in $\Pi$ and take the form
\begin{equation}\label{eq:WsP}
x=C\exp\Big(\dfrac{-\pi k}{\omega}\Big)+o(1)_{k\to +\infty},
\end{equation}
\noindent where $o(1)$ stands for a function of $y$ and $z$ that is uniformly small together with its derivatives up to order $r-2$. Those manifolds accumulate on $\Pi_0$ in $C^0$-topology as $k \to +\infty$. 
\par{}
For any $\mu\neq0$ sufficiently close to 0, there exists a constant $C_1$ such that points in $\{P_k^+\}$ which satisfy the condition
\begin{equation}\label{eq:survivedp}
x_k > C_1|\mu|^{\frac{1}{\rho}} 
\end{equation} 
remain index-1 fixed points of the map $T_1$, and their local stable manifolds $W^s_{loc}(P_k^+)$ take the same form as given by formula (\ref{eq:WsP}).
\end{lem}
%
%
%
\noindent \textit{Proof.} Let $\mu=0$. We first find the fixed points $P^+_k$, which can be done by plugging $(\bar{x}=x,\bar{y}=y,\bar{z}=z)$ into (\ref{eq:1setting_4}). From the last two equations in (\ref{eq:1setting_4}), the coordinates $y$ and $z$ can be expressed as functions of $x$, which leads to the equation for coordinate $x$:
\begin{equation}\label{eq:WsP_1}
x=Ax^\rho\cos(\omega\ln\dfrac{1}{x}+\theta)+o(x^\rho).
\end{equation}
We have the fixed points $P^+_k$ with
\begin{equation}\label{eq:WsP_2}
\begin{array}{rcl}
x_k&=&C\exp\Big(\dfrac{-\pi k}{\omega}\Big)+o\Big(\exp\Big(\dfrac{-\pi k}{\omega}\Big)\Big) ,\\[10pt]
y_k&=&1+o\Big(\exp\Big(\dfrac{-\pi k}{\omega}\Big)\Big) ,\\[10pt]
z_k&=&z^+ +o\Big(\exp\Big(\dfrac{-\pi k}{\omega}\Big)\Big) ,
\end{array}
\end{equation}
\noindent where $y_k,x_k,z_k$ are the coordinates of $P^+_k$, $C=\exp(({2\theta-\pi})/{2\omega})$, and $k$ is any positive integer greater than some sufficiently large $K$. Let us show that the points $P^+_k$ are of index 1. Recall the transformation \eqref{eq:mr1.11} of $x$-coordinate of points on $\Pi_1$:
\begin{equation*}
\xi=\omega\ln\dfrac{1}{x}+\theta-2\pi j \quad \xi\in[0,2\pi)
\end{equation*}
\noindent by which we divide the cross-section into different regions $V_j$ and let $\xi$ be a new coordinate in each region. By Lemma \ref{lem:index_2}, a fixed point $P(x,y,z)$ of $T_i$ $(i=1,2)$ has index 2 only if $\cos\xi$ is close to a value bounded away from zero. However, the first equation in (\ref{eq:WsP_2}) implies that $\cos\xi_k$ is small when $k$ is sufficiently large. We also note that, under our consideration, the index of a periodic point is at most 2 since the multipliers corresponding to $z$ coordinates stay inside the unit circle for all the small perturbations. 
\par{}
We now consider the inverse image under $T_1$ of a small piece of the surface $\{x=x_k\}$ containing $P^+_k$. By formula (\ref{eq:1setting_4}), we have
\begin{equation}\label{eq:WsP_3}
\sin\Big(\dfrac{\pi}{2}-\theta-\omega\ln\dfrac{1}{x}\Big)=\dfrac{1}{yA}\Big(\dfrac{x_k}{x^\rho}+o(1)_{x\to 0}\Big) ,
\end{equation}
\noindent where $x$ and $y$ are coordinates of the points in the inverse image ($z$ coordinates are in the $o(1)$ term). Note that $x$ is bounded since the small cross-section $\Pi$ is bounded. We have following equation if $x$ and ${x_k}/{x^\rho}$ are sufficiently small:
\begin{equation}\label{eq:WsP_4}
\dfrac{\pi}{2}-\theta-\omega\ln\dfrac{1}{x}=\dfrac{1}{yA}\Big(\dfrac{x_k}{x^\rho}+o(1)_{x\to 0}\Big)+m\pi\quad ,m=0,\pm1,\pm2,\dots \quad,
\end{equation}
\noindent which, by noting that the surface contains $P^+_k$, leads to
\begin{equation}\label{eq:WsP_5}
x=C\exp\Big(\dfrac{-\pi k}{\omega}\Big)+o(1)_{k\to +\infty} ,
\end{equation}
\noindent where the term $o(1)$ stands for a function of $y$ and $z$ that is uniformly small together with its derivatives up to order $r-2$. Formula (\ref{eq:WsP_5}) is valid for all values of $y,z$, where $(x,y,z)\in\Pi$, if $x$ and ${x_k}/{x^\rho}$ are sufficiently small. This requirement is equivalent to that $k$ is sufficiently large. One can check that the successive backward iterates of a small piece of the surface $x=x_k$ containing $P^+_k$ take the same form as (\ref{eq:WsP_5}), where the term $o(1)$ stays uniformly small. Since $W^s_{loc}(P)$ is the limit of a sequence of those iterates, $W^s_{loc}(P)$ is given by (\ref{eq:WsP_5}). Obviously, those manifolds accumulate on $\Pi_0$ as $k \to +\infty$ in $C^0$-topology. 
%
%
\par{} At $\mu\neq 0$, the fixed points are given by 
\begin{equation}\label{eq:WsP_6}
x=\mu+Ax^\rho\cos(\omega\ln\dfrac{1}{x}+\theta)+o(x^\rho),
\end{equation}
\noindent which will still lead to the formula (\ref{eq:WsP_2}) for the coordinates of the fixed point $P^+_k$ if $|\mu|\exp({\pi\rho k}/{\omega})$ is sufficiently small, i.e. there is a sufficiently large constant $C_1$ such that $x_k < C_1|\mu|^{\frac{1}{\rho}}$. Note that formula (\ref{eq:WsP_2}) now only gives us finitely many fixed points and they are still of index 1.
\par{}
We now consider again the inverse image under $T_1$ of a small piece of the surface $\{x=x_k\}$ containing $P^+_k$ but with $\mu\neq 0$. From equations in \eqref{eq:1setting_4}, we have
\begin{equation}\label{eq:WsP_7}
\sin\Big(\dfrac{\pi}{2}-\theta-\omega\ln\dfrac{1}{x}\Big)=\dfrac{1}{yA}\Big(\dfrac{x_k}{x^\rho}-\dfrac{\mu}{x^\rho}+o(1)_{x\to 0}\Big) ,
\end{equation}
\noindent where $(x,y)$ are coordinates of the points in the inverse image and $y$ is bounded since the small cross-section $\Pi$ is bounded. In addition to $x$ and ${x_k}/{x^\rho}$, if ${|\mu|}/{x^\rho}$ is also sufficiently small, then we have the following equation similar to equation (\ref{eq:WsP_4}): 
\begin{equation}\label{eq:WsP_8}
\dfrac{\pi}{2}-\theta-\omega\ln\dfrac{1}{x}=\dfrac{1}{yA}\Big(\dfrac{x_k}{x^\rho}-\dfrac{\mu}{x^\rho}+o(1)_{x\to 0}\Big)+m\pi\quad ,m=0,\pm1,\pm2,\dots \quad,
\end{equation}
\noindent which also gives
\begin{equation}\label{eq:WsP_9}
x=C\exp\Big(\dfrac{-\pi k}{\omega}\Big)+o(1)_{k\to +\infty} .
\end{equation}
\noindent Formula (\ref{eq:WsP_9}) has the same form as (\ref{eq:WsP_5}), and it is valid for all values of $y,z$, where $(x,y,z)\in\Pi$, if $|\mu|\exp({\pi\rho k}/{\omega}),x,{x_k}/{x^\rho}$ and ${\mu}/{x^\rho}$ are sufficiently small. This is equivalent to that $k$ is sufficiently large and $|\mu|\exp({\pi\rho k}/{\omega})$ is sufficiently small. This can be achieved since $\rho<1$ and we can choose sufficiently small $\mu$ and sufficiently large $k$ independently. It can be checked that the successive backward iterates of the curve given by (\ref{eq:WsP_9}) take the same form, where the term $o(1)$ stays uniformly small. It follows that the local stable manifold $W^s_{loc}(P)$ is given by (\ref{eq:WsP_9}).\qed
\par{}
For a non-zero $\mu$ value, we can no longer chose $P$ as close to the surface $\Pi_0$ as we want, since infinitely many points in $\{P^+_k\}$ which accumulate on $\Pi_0$ are destroyed. In order to find the intersection $W^s(P)\cap W^u(Q)$, more details are required on the positions of the points $P_k$ that remain index-1 and fixed after changing $\mu$.
\par{}
Let $P^+_{k^*}\in\{P^+_k\}$ be the first point (the one with the largest subscript) in $\{P^+_k\}$ satisfying $x_{k} > C_1|\mu|^{\frac{1}{\rho}}$ where $\mu\neq 0$ (see Lemma \ref{lem:WsP}). Obviously, the number $k^*$ mainly depends on $\mu$, and $k^*\to +\infty$ as $\mu \to 0$. We have the following results.
\begin{lem}\label{lem:ubdyk*}
For any given constant $C_2>0$, there exists a positive number $\mu(C_2)$ such that the inequality 
\begin{equation}
x_{k^*} < C_2|\mu|^{\frac{1}{2\rho}}
\end{equation}
\noindent holds for all $|\mu|<\mu(C_2)$.
\end{lem}
\noindent\textit{Proof.}
By the definition of the point $P^+_{k^*}$, it is sufficient to prove that, for any given $C_2$, there exists a point $P^+_k$ satisfying
\begin{equation}\label{eq:upbyk*}
C_1|\mu|^{\frac{1}{\rho}}<x_k<C_2|\mu|^{\frac{1}{2\rho}},
\end{equation}
\noindent where $C_1$ is the constant in Lemma \ref{lem:WsP}. Recall that we have 
\begin{equation*}
x_k=C\exp\Big(\dfrac{-\pi k}{\omega}\Big)+o\Big(\exp\Big(\dfrac{-\pi k}{\omega}\Big)\Big),
\end{equation*} 
\noindent from formula (\ref{eq:WsP_2}). Now by letting 
\begin{equation}\label{eq:upbyk*_1}
\mu=\exp\bigg(\dfrac{-\pi j-\xi_\mu+\eta}{\omega}\bigg)
\end{equation}
\noindent we have 
\begin{equation*}
\frac{x_k}{|\mu|^{\frac{1}{\rho}}}=\left|\frac{x_k^\rho}{\mu}\right|^{\frac{1}{\rho}}=
C\exp\bigg(\dfrac{\xi_\mu-\eta}{\omega}\bigg)^{\frac{1}{\rho}}\exp\Bigg(\dfrac{(j-\rho k)\pi}{\omega}\Bigg)^{\frac{1}{\rho}}+\dots
\end{equation*} 
\noindent and
\begin{equation*}
\dfrac{|\mu|^{\frac{1}{2\rho}}}{x_k}=\left|\frac{\mu}{x_k^{2\rho}}\right|^{\frac{1}{2\rho}}=
\dfrac{1}{C}\exp\bigg(\dfrac{\eta-\xi_\mu}{\omega}\bigg)^{\frac{1}{2\rho}}\exp\Bigg(\dfrac{(2\rho k-j)\pi}{\omega}\Bigg)^{\frac{1}{2\rho}}+\dots,
\end{equation*} 
\noindent where dots denote small terms that go to zero as $j$ and $k$ go to plus infinity. Consequently, there exists sufficiently large integer $K$ and $J$ such that the inequalities
\begin{equation*}
\frac{x_k}{|\mu|^{\frac{1}{\rho}}}>C'_1 \exp\Bigg(\dfrac{(j-\rho k)\pi}{\omega}\Bigg)^{\frac{1}{\rho}}
\quad \mbox{and} \quad
\dfrac{|\mu|^{\frac{1}{2\rho}}}{x_k}>
C'_2 \exp\Bigg(\dfrac{(2\rho k-j)\pi}{\omega}\Bigg)^{\frac{1}{2\rho}},
\end{equation*}
\noindent hold for all $j>J$ and $k>K$, where $C'_1$ and $C'_2$ are two constants that do not depend on $j$ and $k$. It follows that, in order to obtain inequality (\ref{eq:upbyk*}), it is now sufficient to find $j$ and $k$ that satisfy inequalities
\begin{equation*}
j-\rho k>\dfrac{\omega\rho\ln C_1/C'_1}{\pi} \quad \mbox{and} \quad 2\rho k-j>-\dfrac{2\omega \rho\ln C_2C'_2}{\pi} ,
\end{equation*}
\noindent i.e.
\begin{equation}\label{eq:upbyk*_2}
\dfrac{\omega\rho\ln C_1/C'_1}{\pi} + \rho k<j<\dfrac{2\omega \rho\ln C_2C'_2}{\pi} +2\rho k .
\end{equation}
\noindent Obviously, for any sufficiently large integer $j$, one can find an integer $k$ such that the above inequality holds (note $2\rho<1$).
\par{}
Let $j_0$ be such that for every $j>j_0$ there exists $k$ such that $(j,k)$ is a solution to inequality \eqref{eq:upbyk*_2}. Now let
\begin{equation}
\mu(C_2)=\exp(\dfrac{-\pi j_0+\theta}{\omega}).
\end{equation} 
\noindent Then, for any $|\mu|<\mu(C_2)$, the corresponding $j$ given by equation \eqref{eq:upbyk*_1} satisfies $j\geqslant j_0$. The lemma is proven.\qed
\begin{lem}\label{lem:homorelated}
There exists a constant $K$ for all $\rho$ values close to $\rho^*$ (the $\rho$ value of the original system) such that points in $\{P^+_k\}$ with $K(\rho)<k\leq k^*$ are homoclinically related.
\end{lem}
\noindent\textit{Proof.} The case where $\mu=0$ (i.e. $k^*=+\infty$) is the result of \citep{sh70}. When $\mu$ is non-zero, as long as $y_k < C_1|\mu|^{\frac{1}{\rho}}$, the parameter $\mu$ will enter small terms of all equations used in the computation in \citep{sh70}, and therefore we have Lemma \ref{lem:homorelated}.\qed
\subsubsection{The transverse intersection}\label{prf1.5.2}
If the intersection $W^s(P)\cap W^u(Q)$ exists, then it must be transverse. This is because that the invariant manifold $\cap W^u(Q)$ is transverse to the strong-stable foliation, but $W^s(P)$ consists of leaves of the strong-stable foliation.
\par{}
We assume that the point $P\in\{P^+_k\}$ used in the previous sections has a subscript satisfying
\begin{equation}\label{eq:transsymmetric1}
K<k<k^*,
\end{equation}
where $k^*$ is the largest subscript among points in $\{P^+_k\}$ which remain index-1 and fixed after a small perturbation, and $K$ is the constant given by Lemma \ref{lem:homorelated}. This means that our point $P$ lies in the set of the remaining points in $\{P^+_k\}$ which are not only index-1 and fixed but also homoclinically related with each other after small perturbations in $\rho$. 
\par{}
In what follows we will show that, for any periodic point $Q$ of $T$ with index 2, its unstable manifold $W^u(Q)$ transversely intersects the stable manifold $W^s(P)$. To achieve this, we use the following result:
\begin{lem}\label{lem:WuQ}
The unstable manifold $W^u(Q)$ of the orbit of an index-2 periodic point $Q$ intersects $\Pi_0=\{y=0\}$ transversely.
\end{lem}
This is Lemma 5 in \citep{li16} and here we only sketch the proof. It is proved by considering the quotient map $\tilde{T}\equiv(\tilde{T}_1,\tilde{T}_2)$ obtained from the Poincaré map $T\equiv(T_1,T_2)$ by taking quotient along leaves of the strong-stable foliation $\mathcal{F}_1$ on $\Pi$. The map $\tilde{T}$ acts on the 2-dimensional quotient cross-section $\tilde{\Pi}=\Pi\cap \{z=z^*\}$ where $z^*$ is some constant and $\|z^*\|<\delta$. For any region $V\subset\tilde{\Pi}$, its image $\tilde{T}(V)$ is the projection of $T(V)$ onto $\tilde{\Pi}$ along the leaves of $\mathcal{F}_1$. Let $\tilde{Q}$ be the projection of $Q$ on $\tilde{\Pi}$ by the leaf which goes through $Q$. Note that $W^u(\tilde{Q})$ is obtained by taking limit of the iterates of a small two-dimensional neighbourhood of $\tilde{Q}$ on $\tilde{\Pi}$. By the volume-hyperbolicity of the flow and the absolute continuity of the foliation, two-dimensional areas on $\tilde{\Pi}$ is expanding under $\tilde{T}$. Therefore, the unstable manifold $W^u(\tilde{Q})$ intersects $W^s(O)\cap\tilde{\Pi}$, which implies the lemma. 
\par{}
By Lemma \ref{lem:WuQ}, we can take a connected component $L\subset W^u_{loc}(Q)$ such that it intersects $\Pi_0$. We will consider the iterate of $L$ under the map $T^{(2)}$. We show (see Lemma \ref{lem:WuQ_2} below) that some iterate of $L$ transversely intersects the local stable manifold $W^s_{loc}(P^+_{k^*})$ if $\Gamma^+$ and $\Gamma^-$ form double-round homoclinic loops, which is the case when we consider the parameters given by Lemma \ref{lem:quasisymmetric}. 
\par{}
Recall that the point $P$ picked is homoclinically related to $P^+_{k^*}$. Therefore, by $\lambda$-lemma, we obtain the transverse intersection $W^s(P)\cap W^u(Q)$. This intersection along with the quasi-transverse intersection $W^u(P)\cap W^s(Q)$ obtained in Lemma \ref{lem:quasisymmetric} immediately implies the existence of a heterodimensional cycle of the Poincaré map $T$ associated to $P$ and $Q$. It follows that we obtain a heterodimensional cycle associated to two periodic orbits of indices 2 and 3 in the full system $X_{\mu_j^n,\rho_j^n,\nu_j^n}$. In this way, Theorem \ref{thm:hetero_1} will be proven. 
\par{}
Let us now give the lemma used in the above argument. Let $l\subset L$ be a curve joining two points $M_0(0,x_0,z_0)\in L\cap\Pi_0$ and $M_1(x_1,y_1,z_1)\in (L\cap\Pi_1)\setminus\Pi_0$. We now consider the iterate of this curve $l$ under the map $T$. 
\begin{lem}\label{lem:WuQ_2}
If $\Gamma^+$ and $\Gamma^-$ are two double-round homoclinic loops, then there exists some $i$ such that the iterate $T^{(i)}(l)$ intersects the local stable manifold $W^s_{loc}(P^+_{k^*})$. 
\end{lem}
\noindent This lemma immediately implies the existence of the non-empty transverse intersection $W^s_{loc}(P^+_{k^*})\cap W^u_{loc}(Q)$.
%
%
\par{}
\noindent\textit{Proof of Lemma \ref{lem:WuQ_2}.} 
By Lemma \ref{lem:WsP}, the local stale manifold $W^s_{loc}(P^+_{k^*})$ is given by
\begin{equation}\label{eq:WuQ_2_1}
x=C\exp(\dfrac{-\pi k^*}{\omega})+o(1)_{k^*\to +\infty}=x_{k^*}+o(1)_{k^*\to +\infty},
\end{equation}
\noindent where $x_{k^*}$ is the $x$-coordinate of $P^+_{k^*}$. Therefore, by Lemma \ref{lem:ubdyk*}, we can find two positive constants $C_2$ and $\mu(C_2)$ such that $W^s_{loc}(P^+_{k^*})$ is below the surface $\{x=C_2|\mu|^\frac{1}{2\rho}\}$ for all $\mu\in(-\mu(C_2),\mu(C_2))$. Consequently, if $x_1\geqslant C_2|\mu|^\frac{1}{2\rho}$, then Lemma \ref{lem:WuQ_2} automatically holds. We now assume 
\begin{equation}\label{eq:WuQ_2_assumption}
x_1< C_2|\mu|^\frac{1}{2\rho}.
\end{equation}
\par{}
We first show that there exists a point $M_2(y_2,x_2,z_2)\in l$ with $0<x_2<x_1$ such that the $x$-coordinate $\bar{\bar{x}}_2$ of its second iterate $T^{(2)}(M_2)$ is lager than $x_1$. Note that we have $T^{(2)}(M_2)=T_1^{(2)}(M_2)$ if $\mu>0$ and $T^{(2)}(M_2)=T_2\circ T_1(M_2)$ if $\mu<0$. For certainty, we consider the case where $\mu>0$. The same result holds for the other case.
\par{}
Since $M_0\in\Pi_0$ and the homoclinic loop $\Gamma^+$ is double-round (i.e. we consider here only the parameter values given by Lemma \ref{lem:quasisymmetric}), we have 
\begin{equation}\label{eq:WuQ_2_2}
T_1(M_0)=M^+(\mu,1,z^+) \quad \mbox{and} \quad T_1^{(2)}(M_0)=(0,\bar{\bar{y}}_0,\bar{\bar{z}}_0).
\end{equation} 
\noindent Let $M(x,y,z)$ be an arbitrary point on $l$ and we consider the $x$-coordinate $\bar{\bar{x}}$ of its second iterate $T_1^{(2)}(M)$.  
Recall the equation for $x$-coordinate in the formula (\ref{eq:1setting_4}) of the map $T_1$, which is 
\begin{equation}\label{eq:WuQ_2_3}
F(x,y,z):=\bar{x}=\mu+Ayx^\rho\cos{(\omega \ln {\dfrac{1}{x}}+\theta)}+o(x^\rho).
\end{equation}
\noindent 
By mean value theorem, we have 
\begin{equation}\label{eq:WuQ_2_4}
\begin{array}{rcl}
\bar{\bar{x}}&=&\bar{\bar{x}}-0\\[10pt]
&=&F(\bar{x},\bar{y},\bar{z})-F(\mu,1,z^+) \\[10pt]
&=&\dfrac{\partial F(x_t,y_t,z_t)}{\partial x}(\bar{x}-\mu)+\dfrac{\partial F(x_t,y_t,z_t)}{\partial y}(\bar{y}-1)
+\dfrac{\partial F(x_t,y_t,z_t)}{\partial z}(\bar{z}-z^+),
\end{array}
\end{equation}
\noindent where $(x_t,y_t,z_t)=(1-t)(\bar{x},\bar{y},\bar{z})+t(\mu,1,z^+)$ for some $t\in(0,1).$ By equation (\ref{eq:WuQ_2_3}) and formula (\ref{eq:1setting_4}) of the map $T_1$, equation (\ref{eq:WuQ_2_4}) yields
\begin{equation}\label{eq:WuQ_2_5}
\begin{array}{rcl}
\bar{\bar{x}}&=&(\sqrt{\rho^2+\omega^2}x_t^{\rho-1}\cos(\omega\ln\dfrac{1}{x_t}+\theta-\varphi)+o(x_t^{\rho-1}))(\bar{x}-\mu)+O(x_t^\rho)(\bar{x}-1)\\[10pt]
&&+O(x_t^\rho)(\bar{z}-z^+) \\[10pt]
&=&\sqrt{\rho^2+\omega^2}x_t^{\rho-1}\cos(\omega\ln\dfrac{1}{x_t}+\theta-\varphi)Ayx^\rho\cos(\omega\ln\dfrac{1}{x}+\theta)+o(x_t^{\rho-1}x^\rho),
\end{array}
\end{equation}
\noindent where $\phi=\arctan(\dfrac{\omega}{\rho})$.
\par{}
Let us now find an estimate for $\bar{\bar{x}}$. It can be seen from equation (\ref{eq:WuQ_2_5}) that the first term in (\ref{eq:WuQ_2_5}) is dominant if 
\begin{equation}\label{eq:WuQ_2_5.1}
\cos(\omega\ln\dfrac{1}{x}+\theta)\neq 0 \quad\mbox{and}\quad \cos(\omega\ln\dfrac{1}{x_t}+\theta-\varphi)\not\equiv 0  \quad t\in(0,1).
\end{equation}
Moreover, this first term is also monotone if $\cos(\omega\ln{x_t}^{-1}+\theta-\varphi)$ does not change sign for all $t\in(0,1)$. In what follows we find points on the curve $l$ satisfying these conditions and the lower bound of the $x$-coordinates $\bar{\bar{x}}$ of their second iterate under the map $T_1$.
\par{}
Obviously, there exists a sufficiently small $\varepsilon>0$ such that if 
\begin{equation}\label{eq:WuQ_2_6}
\left|\ln\dfrac{1}{\mu}-\ln\dfrac{1}{\bar{x}}\right|<\ln(1+\varepsilon)\quad \mbox{i.e.}\quad 
\left|\dfrac{\bar{x}}{\mu}-1\right|<\varepsilon
,
\end{equation} 
\noindent then $\cos(\omega\ln{x_t}^{-1}+\theta-\varphi)$ does not change sign for all $t\in(0,1)$, and therefore equation (\ref{eq:WuQ_2_5}) implies that 
\begin{equation}\label{eq:WuQ_2_7}
\left|\bar{\bar{x}}\right|\geqslant C_3\mu^{\rho-1}|Ayx^\rho\cos(\omega\ln\dfrac{1}{x}+\theta)+o(x^\rho)|,
\end{equation}
\noindent where $C_3$ is a constant depending on $\varepsilon$. Now recall the variable $\xi$ introduced by the relation (\ref{eq:mr1.11}) that
\begin{equation}\label{eq:WuQ_2_7.1}
\omega\ln\dfrac{1}{|x|}=2\pi j+\xi-\theta \quad \xi \in[0,2\pi).
\end{equation} 
\noindent We consider a sequence $\{M_{cj}(x_{cj},y_{cj},z_{cj})\}_j$ of points on the curve $l$ such that  
\begin{equation}\label{eq:WuQ_2_8}
\cos(\omega\ln\dfrac{1}{x_{cj}}+\theta)=\cos\xi_{cj}\neq 0 .
\end{equation}
\noindent Specifically, we consider points that satisfy
\begin{equation}\label{eq:WuQ_2_9}
\dfrac{\partial F(x_{cj},y_{cj},z_{cj})}{\partial x}=0
\end{equation}
\noindent i.e.
\begin{equation*}
\sqrt{\rho^2+\omega^2}x_{cj}^{\rho-1}\cos(\xi_{cj}-\varphi)+o(x_{cj})=0,
\end{equation*}
\noindent which implies
\begin{equation}\label{eq:WuQ_2_10}
\xi_{cj}=\dfrac{\pi}{2}+\phi+k\pi+\cdots
\end{equation}
\noindent where $k=0,1$ since $\xi_{cj}\in[0,2\pi)$ and the dots denote terms that go to zero as $j$ goes to plus infinity. Generically, inequality (\ref{eq:WuQ_2_8}) is satisfied by those $\xi_{cj}$. 
\par{}
Note that, for any $\varepsilon$, condition (\ref{eq:WuQ_2_6}) can be satisfied by all points in $\{M_{cj}\}$ with sufficiently large $j$. This is because, by equation (\ref{eq:WuQ_2_3}), we have
\begin{equation}\label{eq:WuQ_2_11}
\left|\dfrac{\bar{x}_{cj}}{\mu}-1\right|=\left|\dfrac{Ay_{cj}x_{cj}^\rho\cos\xi_{cj}+o(x_{cj}^\rho)}{\mu}\right|<\varepsilon
\end{equation}
\noindent when $x_{cj}$ is sufficiently small (i.e. $j$ is sufficiently large) since $y_{cj}$ is uniformly bounded by the definition of the cross-section $\Pi$. Indeed, by the relation (\ref{eq:WuQ_2_7.1}), equation (\ref{eq:WuQ_2_11}) yields
\begin{equation}\label{eq:WuQ_2_12}
\left|Ay_{cj}\exp\bigg(\dfrac{-2\pi\rho j-\rho\xi_{cj}+\rho\theta}{\omega}\bigg)\cos\xi_{cj}+o\bigg(\exp\bigg(\dfrac{-2\pi\rho j}{\omega}\bigg)\bigg)\right|<\mu\varepsilon ,
\end{equation} 
\noindent i.e.
\begin{equation}\label{eq:WuQ_2_13}
j>\dfrac{\theta-\xi_{cj}}{2\pi}-\dfrac{\omega}{2\pi\rho}\ln\dfrac{\mu\varepsilon}{Ay_{cj}|\cos\xi_{cj}|}+o(1)_{j\to +\infty}=:J_1+o(1)_{j\to +\infty}.
\end{equation}
\noindent Now let $J_2$ be the smallest integer such that the above small term $o(1)_{j\to +\infty}$ is lesser than $1$ for all $j>J_2$, and let $J=\max(J_1+1,J_2)$. It follows that all points $M_{cj}$ with $j>J$ satisfy the condition (\ref{eq:WuQ_2_6}). 
\par{}
Note that equation (\ref{eq:WuQ_2_8}) implies that there exists a constant $C_4$ which does not depend on $\mu$ such that we have 
\begin{equation}\label{eq:WuQ_2_14}
|Ay_{cj}x_{cj}^\rho\cos\xi_{cj}+o(x_{cj}^\rho)|>C_4 x_{cj}^\rho.
\end{equation}
\noindent Thus, for points $M_{cj}$ with $j>J$, inequality (\ref{eq:WuQ_2_7}) now implies 
\begin{equation}\label{eq:WuQ_2_15}
\left|\bar{\bar{x}}_{cj}\right|> C_3 C_4\mu^{\rho-1}x_{cj}^\rho.
\end{equation}
\par{}
We next claim that there exists a positive constant $C_5$ independent of $\mu$ such that one can always find some $j_0>J$ satisfying the inequality
\begin{equation}\label{eq:WuQ_2_16}
x^\rho_{cj_0}>\dfrac{\varepsilon\mu}{C_5}
\end{equation}
\noindent at $\mu\neq 0$. By the relation (\ref{eq:WuQ_2_7.1}), inequality (\ref{eq:WuQ_2_16}) is equivalent to
\begin{equation}\label{eq:WuQ_2_17}
\exp\bigg(\dfrac{-2\pi\rho j_0-\rho\xi_{cj_0}+\rho\theta}{\omega}\bigg)>\dfrac{\varepsilon\mu}{C_5} ,
\end{equation} 
\noindent i.e.
\begin{equation}\label{eq:WuQ_2_18}
j_0<\dfrac{\theta-\xi_{cj_0}}{2\pi}-\dfrac{\omega}{2\pi\rho}\ln\dfrac{\mu\varepsilon}{C_5}.
\end{equation}
\noindent By comparing inequality (\ref{eq:WuQ_2_18}) with the definition of $J$ (\ref{eq:WuQ_2_13}) and noting that $Ay_{cj}|\cos\xi_{cj}|$ is bounded, one can easily find the constant $C_5$ stated in the claim.
\par{}
We now consider the point $M_{cj_0}$. The assumption $x_1< C_2|\mu|^\frac{1}{2\rho}$ (\ref{eq:WuQ_2_assumption}) and inequality (\ref{eq:WuQ_2_16}) imply
\begin{equation}\label{eq:WuQ_2_19}
x^\rho_{cj_0}>\dfrac{\varepsilon\mu^\frac{1}{2}}{C_2 C_5}x_1^\rho.
\end{equation}
\noindent The $x$-coordinate $\bar{\bar{x}}_{cj_0}$ of the second iterate $T_1^{(2)}(M_{cj_0})$ can be now estimated by inequalities (\ref{eq:WuQ_2_15}) and (\ref{eq:WuQ_2_19}) as 
\begin{equation}\label{eq:WuQ_2_20}
\left|\bar{\bar{x}}_{cj_0}\right|> \dfrac{\varepsilon C_3 C_4}{C_2 C_5} \mu^{\rho-1}x_{1}^\rho.
\end{equation}
\noindent Note that the constants $C_2,C_3,C_4$ and $C_5$ do not depend on $\mu$, and the point $M_{cj_0}$ exists for all $\mu\neq 0$ by the above claim. Therefore, we can choose $\mu$ sufficiently small such that inequality (\ref{eq:WuQ_2_20}) implies
\begin{equation}\label{eq:WuQ_2_21}
\left|\bar{\bar{x}}_{cj_0}\right|> 2 x_1.
\end{equation}
\par{}
Note that the sign of $\bar{\bar{y}}_{cj_0}$ is the same of that of $\cos(\omega\ln{y_t}^{-1}+\theta-\varphi)\cos\xi_{cj_0}$ (see (\ref{eq:WuQ_2_5})). When $\varepsilon$ is small, the sign of $\cos(\omega\ln{y_t}^{-1}+\theta-\varphi)$ depends on the value of $\mu$. Also, it can be seen from equation (\ref{eq:WuQ_2_10}) that 
\begin{equation*}
Ay_{cj_0}x_{cj_0}^\rho\cos\xi_{cj_0}+o(x_{cj_0}^\rho)>0 \quad \mbox{if}\quad k=1,
\end{equation*}
\noindent and
\begin{equation*}
Ay_{cj_0}x_{cj_0}^\rho\cos\xi_{cj_0}+o(x_{cj_0}^\rho)<0 \quad \mbox{if}\quad k=0.
\end{equation*}
\noindent It follows that we can choose $k$ accordingly such that $\bar{\bar{x}}_{cj_0}$ given by equation (\ref{eq:WuQ_2_5}) is positive. Consequently, we can rewrite inequality (\ref{eq:WuQ_2_21}) as
\begin{equation}\label{eq:WuQ_2_22}
\bar{\bar{x}}_{cj_0}>2x_1.
\end{equation}
\par{}
Now Let $M_{cj_0}$ be the point $M_2(x_2,y_2,z_2)$ mentioned in the beginning of the proof. We return to the iterate of the curve $l$ joining points $M_0$ and $M_1(x_1,y_1,z_1)$. From the above argument, the second iterate $T^{(2)}_1(l)$ contains a curve $l_1$ joining points $T^{(2)}_1(M_0)=:M_3$ and $T^{(2)}_1(M_2)=:M_4(x_4,y_4,z_4)$ such that $x_4>2x_1$. Note that we have $M_3\in\Pi_0$ since the separatrices $\Gamma^+$ and $\Gamma^-$ form double-round homoclinic loops for the $\mu$ values considered, i.e. $T^{(2)}(\Pi_0)\subset \Pi_0$. Hence, we can apply the same argument to the curve $l_1$, and obtain a point $M_5\in l_1$ and a new curve $l_2\subset T^{(2)}(l_1)$ joining points $T^{(2)}(M_3)=:M_6\in\Pi_0$ and $T^{(2)}(M_5)=:M_7(x_7,y_7,z_7)$ such that $x_7>2x_4$. This procedure can be continued until we find a curve $l_k$ joining points $M_{3k}\in\Pi_0$ and $M_{3k+1}(x_{3k+1},y_{3k+1},z_{3k+1})$ such that the assumption (\ref{eq:WuQ_2_assumption})
\begin{equation*}
x_{3k+1} < C_2|\mu|^\frac{1}{2\rho}
\end{equation*}
\noindent is violated. Consequently, we now have
\begin{equation*}
x_{3k+1}\geqslant C_2|\mu|^\frac{1}{2\rho},
\end{equation*}
which means that $l_k$ intersects $W^s_{loc}(P^+_{k^*})$. The lemma is proved. \qed
\subsection{Proof of Lemma \ref{lem:WsQ}}\label{sec:prf1.6}
\noindent{\it Proof.} Let $M_1(x_1,y_1,z_1)$ be a point on $\Pi_1$ and denote by $M_2(x_2,y_2,z_2)$ its image under $T_1$. We first show that there exist two cones $\mathcal{C}_1$ and $\mathcal{C}_2$ at $M_1$ and $M_2$ respectively such that the preimage of any tangent vector in $\mathcal{C}_2$ under $T_1$ lies in $\mathcal{C}_1$, provided $y_1$ is sufficiently small.
\par{}
We start by estimating the norms of the preimages. By formula \eqref{eq:WsQ_1} for the derivative of $T_1$, we obtain
\begin{equation}\label{eq:WsQ_1_new}
\begin{array}{c}
\left.\mbox{D}T_1\right|_{M_1}=\\[15pt]
\left(
\begin{array}{lll}
Ax_1 ^{\rho-1}(\rho\cos\xi_1+\omega\sin\xi_1)+o(x_1 ^{\rho-1}) & Ax_1 ^\rho\cos\xi_1+o(x_1 ^\rho) &\boldsymbol{a} \\[10pt]

\begin{array}{l}
-A_1x_1 ^{\rho -1}(\rho \cos(\xi_1+\eta_1-\eta)\\
+\omega\sin(\xi_1+\eta_1-\eta))+o(x_1 ^{\rho -1})
\end{array}
&
\begin{array}{l}
A_1x_1 ^\rho\cos(\xi_1+\eta_1-\eta)\\
+o(x_1 ^\rho)
\end{array}
 & \boldsymbol{a}_1\\[20pt]

\begin{array}{l}
-A_2x_1 ^{\rho -1}(\rho \cos(\xi_1+\eta_2-\eta)\\
+\omega\sin(\xi_1+\eta_2-\eta))+o(x_1 ^{\rho -1})
\end{array}
&
\begin{array}{l}
A_2x_1 ^\rho\cos(\xi_1+\eta_2-\eta)\\
+o(x_1 ^\rho)
\end{array}
 &\boldsymbol{a}_2\\[10pt]

\dots&\dots&\dots \\[10pt]

\begin{array}{l}
-A_{n-2}x_1 ^{\rho -1}(\rho \cos(\xi_1+\eta_{n-2}-\eta)\\
+\omega\sin(\xi_1+\eta_{n-2}-\eta))+o(x_1 ^{\rho -1})
\end{array}
&
\begin{array}{l}
A_{n-2}x_1 ^\rho\cos(\xi_1+\eta_{n-2}-\eta)\\
+o(x_1 ^\rho)
\end{array}
 &\boldsymbol{a}_{n-2}
 \\

\end{array}
\right),\\
\end{array}
\end{equation}\\[5pt]
\noindent where $\boldsymbol{a}$ and $\boldsymbol{a}_i$ are $1\times (n-3)$ vectors of the form $(o(x_1^\rho),\dots,o(x_1^\rho))$. We rewrite the above matrix as\\[-5pt]
\begin{equation*}
\begin{pmatrix}
y_1^{\rho-1}b_{11} & y_1^{\rho}b_{12} & b_{13} \\
y_1^{\rho-1}b_{21} & y_1^{\rho}b_{22} & b_{23} \\
y_1^{\rho-1}b_{31} & y_1^{\rho}b_{32} & b_{33} 
\end{pmatrix},
\end{equation*}\\[-5pt]
\noindent where $b_{31}$ and $b_{32}$ are $(n-3)\times 1$ vectors, $b_{13}$ and $b_{23}$ are $1\times (n-3)$ vectors, $b_{33}$ is a $(n-3)\times(n-3)$ matrix, and $b_{ij}$ are uniformly bounded.
\par{} Denote by $E$ the $2\times2$ block in the top-left corner of $\left.\mbox{D}T_1\right|_{M_1}$. One can check that the determinant of $E$ is
\begin{equation}
\omega AA_1x_1^{2\rho -1}\sin\eta+o(x_1^{2\rho -1}).
\end{equation}
\noindent We then have
\begin{equation}\label{eq:WsQ_2}
\begin{array}{rcl}
E^{-1}&=&\dfrac{1}{\mbox{det}E}
\left(\begin{array}{ll}
A_1x_1^{\rho}\cos(\xi_1+\eta_1-\eta)+o(x_1^\rho) & 
Ay_1^{\rho}\cos\xi_1+o(y_1^\rho) \\[20pt]

\begin{array}{l}
A_1x_1^{\rho -1}(\rho \cos(\xi_1+\eta_1-\eta)\\
+\omega\sin(\xi_1+\eta_1-\eta))+o(x_1^{\rho -1})
\end{array}
 &
 \begin{array}{l}
 Ax_1^{\rho -1}(\rho \cos\xi_1\\
 +\omega\sin\xi_1)+o(x_1^{\rho-1}))
 \end{array}

\end{array}\right)\\[50pt]
&=&
\begin{pmatrix}
a_{11}y_1^{1-\rho} & a_{12}y_1^{1-\rho} \\
a_{21}y_1^{-\rho} & a_{22}y_1^{-\rho}
\end{pmatrix} ,
\end{array} 
\end{equation}
\noindent where $a_{ij}=\tilde{a}_{ij}+o(1)_{y_1 \to 0}$ while $\tilde{a}_{ij}$ are uniformly bounded when $x_1$ is small. Let $(\Delta x_2,\Delta y_2,\Delta z_2)$ be a vector in the cone $\mathcal{C}_2$, i.e. $|\Delta x_2,\Delta y_2|\leqslant K\|\Delta z\|$ for some given $K>0$. We have\\[-5pt]
\begin{equation}\label{eq:WsQ_3}
\begin{pmatrix}
\Delta x_2 \\
\Delta y_2 \\
\Delta z_2 
\end{pmatrix}
=
\left.\mbox{D}T_1\right|_{M_1}
\begin{pmatrix}
\Delta x_1 \\
\Delta y_1 \\
\Delta z_1 
\end{pmatrix} ,
\end{equation}  \\[-5pt]
\noindent which implies\\[-5pt]
\begin{equation}\label{eq:WsQ_4}
\begin{pmatrix}
\Delta x_1 \\
\Delta y_1 
\end{pmatrix}
=
E^{-1}\left(\begin{pmatrix}
\Delta x_2 \\
\Delta y_2
\end{pmatrix}
-
\begin{pmatrix}
o(x_1^\rho)\Delta z_1 \\
o(x_1^\rho)\Delta z_1
\end{pmatrix}\right).
\end{equation}\\[-5pt]
\noindent After the transformation $\Delta x_1=x_1\Delta u$, the above equation yields\\
\begin{equation}\label{eq:WsQ_5}
\begin{pmatrix}
\Delta u \\
\Delta y_1 
\end{pmatrix}
=
\begin{pmatrix}
a_{11}x_1^{-\rho} & a_{12}x_1^{-\rho} \\
a_{21}x_1^{-\rho} & a_{22}x_1^{-\rho}
\end{pmatrix}
\left(\begin{pmatrix}
\Delta x_2 \\
\Delta y_2
\end{pmatrix}
-
\begin{pmatrix}
o(x_1^\rho)\Delta z_1 \\
o(x_1^\rho)\Delta z_1
\end{pmatrix}\right).
\end{equation}
\par{} We now recover some relations among $\Delta x_1$, $\Delta y_1$ and $\Delta z_1$ from the above equalities. 
By equations (\ref{eq:WsQ_1_new}) and (\ref{eq:WsQ_5}), we get
\begin{equation}\label{eq:WsQ_6}
\begin{array}{rcl}
\Delta z_2 &=& b_{31}x_1^\rho \Delta u+b_{32}x_1^\rho\Delta y_1+o(x_1^\rho)\Delta z_1 \\[10pt]

&=&b_{31}x_1^\rho(a_{11}x_1^{-\rho}\Delta x_2-a_{11}o(1)\Delta z_1+a_{12}x_1^{-\rho}\Delta y_2-a_{12}o(1)\Delta z_1)\\[10pt]

&&+ b_{32}x_1^\rho(a_{21}x_1^{-\rho}\Delta x_2-a_{21}o(1)\Delta z_1+a_{22}x_1^{-\rho}\Delta y_2-a_{22}o(1)\Delta z_1)\\[10pt]&&+o(x_1^\rho)\Delta z_1\\[10pt]

&=&(b_{31}a_{11}+b_{32}a_{21})\Delta x_2+(b_{31}a_{12}+b_{32}a_{22})\Delta y_2+o(x_1^\rho)\Delta z_1,
\end{array}
\end{equation}
\noindent which, by noting that $|\Delta x_2,\Delta y_2|\leqslant K\|\Delta z_2\|$, leads to
\begin{equation}\label{eq:WsQ_7}
\|\Delta z_2\|\leqslant \dfrac{o(x_1^\rho)\|\Delta z_1\|}{1-(\|b_{31}a_{11}+b_{32}a_{21}\|+\|b_{31}a_{12}+b_{32}a_{22}\|)K}.
\end{equation}
\noindent The above inequality along with equation (\ref{eq:WsQ_5}) and the assumption $|\Delta x_2,\Delta y_2|\leqslant K\|\Delta z\|$ implies
\begin{equation}\label{eq:WsQ_8}
\begin{array}{l}
|\Delta x_1|  \leqslant
\left(\left|\dfrac{o(1)(|a_{11}|+|a_{12}|)K}{1-(\|b_{31}a_{11}+b_{32}a_{21}\|+\|b_{31}a_{12}+b_{32}a_{22}\|)K}\right|+o(1)|a_{11}+a_{12}|\right)x_1\|\Delta z_1\|   \\[15pt]

\qquad\,\,\,=:K_1\|\Delta z_1\| 
\end{array}
\end{equation}
\noindent and
\begin{equation}\label{eq:WsQ_9}
\begin{array}{l}
|\Delta y_1|  \leqslant
\left(\left|\dfrac{o(1)(|a_{21}|+|a_{22}|)K}{1-(\|b_{31}a_{11}+b_{32}a_{21}\|+\|b_{31}a_{12}+b_{32}a_{22}\|)K}\right|+o(1)|a_{21}+a_{22}|\right)\|\Delta z_1\|   \\[15pt]

\qquad\,\,\, =:K_2\|\Delta z_1\| ,
\end{array}
\end{equation}
\noindent where $K_1=o(x_1)$, $K_2=o(1)_{x_1\to 0}$ and $K_1,K_2 \leqslant K$ when $x_1$ is sufficiently small. This shows the existence of the desired cones $\mathcal{C}_1$ and $\mathcal{C}_2$ defined in the beginning.
\par{}
Note that the matrices used in the computation above keep the same form if we choose $M_1$ from $\Pi_2$, and therefore all above results hold. This means that for any point $M\in\Pi$, we have a sequence $\{\mathcal{C}_i\}$ of cones along its orbit $\{M_i\}$ such that, for each vector $w\in\mathcal{C}_{i+1}$, its preimage $\mbox{D}T^{-1}(w)$ belongs to $\mathcal{C}_i$, provided $\{M_i\}$ is sufficiently close to $\Pi_0$. In what follows we continue assuming $M_1\in\Pi_1$ and seek for the formula of its local strong-stable manifold.
\par{}
We obtained the above three inequalities (\ref{eq:WsQ_7}) - (\ref{eq:WsQ_9}) by only using the assumption $|\Delta x_2,\Delta y_2|\leqslant K\|\Delta z_2\|$. However, we also know $|\Delta x_2|=k_1\|\Delta z_2\|$ and $|\Delta y_2|=k_2\|\Delta z_2\|$ for some positive constants $k_1$ and $k_2$. By taking into account this fact and assuming $(\Delta x_2,\Delta y_2,\Delta z_2)\in\mathcal{C}_2$, the above computation will lead to the following equalities:
\begin{equation}\label{eq:WsQ_10}
|\Delta x_1|=o(x_1)\|\Delta z_1\| \quad \mbox{and} \quad |\Delta y_1|=o(1)_{x_1\to 0}\|\Delta z_1\|
\end{equation}
\noindent (where the small terms $o(x_1)$ and $o(1)$ are different from those in (\ref{eq:WsQ_8}) and (\ref{eq:WsQ_9})). 
\par{}
Recall the discussion on the non-degeneracy condition in Section \ref{sec:1.1}. We know that there exists a strong-stable foliation $\mathbb{F}_1$ on $\Pi$, and $W^{ss}(M_i)\,(i=1,2)$ are smooth leaves with the form
\begin{equation}\label{eq:WsQ_11}
h_i(z)=
\begin{pmatrix}
h_{i1}(z) \\
h_{i2}(z) \\
z
\end{pmatrix},
\end{equation}
\noindent
where $h_{i1}(z)$ is the $x$-coordinate and $h_{i2}(z)$ is the $y$-coordinate. We also have that ${\mbox{d}h_{ij}}/{\mbox{d} z}$ $(i,j=1,2)$ are uniformly bounded. Let $(h_{11}(z),h_{12}(z),z)$ be a point on $W^{ss}(M_1)$. Denote by $\bar{z}=f_3(h_{11},h_{12},z)=:F(z)$ the third equation in the formula (\ref{eq:1setting_4}) for a point on $W^{ss}(M_1)$, which is the equation for $z$-coordinates in $T_1$.
We take derivative of both sides of $T_1(h_{11}(z),h_{12}(z),z)=(h_{21}(F(z)),h_{22}(F(z)),F(z))$, and obtain
\begin{equation}\label{eq:WsQ_12}
\left.\mbox{D}T_2\right|_{(h_{11}(z),h_{12}(z),z)}
\begin{pmatrix}
\dfrac{\mbox{d} h_{11}(z)}{\mbox{d}z} \\[10pt]
\dfrac{\mbox{d} h_{12}(z)}{\mbox{d}z} \\[10pt]
1
\end{pmatrix}
=
\begin{pmatrix}
\dfrac{\mbox{d} h_{21}(z)}{\mbox{d}z}F' \\[10pt]
\dfrac{\mbox{d} h_{22}(z)}{\mbox{d}z}F' \\[10pt]
F'
\end{pmatrix}.
\end{equation}
\noindent By noting that the derivative ${\mbox{d}h_{2j}}/{\mbox{d} z}$ is uniformly bounded, say $\left\|({\mbox{d}h_{2j}}/{\mbox{d} z})\right\| \leqslant k$, we have $\left\|({\mbox{d} h_{2j}(z)}/{\mbox{d}z})F'\right\|$ $\leqslant k\|F'\|$, which implies $(({\mbox{d} h_{21}(z)}/{\mbox{d}z})F',({\mbox{d} h_{22}(z)}/{\mbox{d}z})F',F')$ is in an above-mentioned cone. Therefore we obtain the following estimate:
\begin{equation}\label{eq:WsQ_13}
\begin{pmatrix}
\dfrac{\mbox{d} h_{11}(z)}{\mbox{d}z} \\[10pt]
\dfrac{\mbox{d} h_{12}(z)}{\mbox{d}z} \\[10pt]
1
\end{pmatrix}=
\begin{pmatrix}
o(h_{11}(z)) \\[10pt]
o(1)_{h_{11}(z) \to 0} \\[10pt]
1
\end{pmatrix}.
\end{equation}
\noindent Moreover, we have
\begin{equation}\label{eq:WsQ_14}
h_1(z)=
\begin{pmatrix}
h_{11}(z)\\
h_{12}(z)\\
z
\end{pmatrix}
=
\begin{pmatrix}
h_{11}(z_1)+\dfrac{\mbox{d} h_{11}(z')}{\mbox{d}z}(z-z_1) \\[10pt]
h_{12}(z_1)+\dfrac{\mbox{d} h_{12}(z'')}{\mbox{d}z} (z-z_1)\\[10pt]
z
\end{pmatrix}=:
\begin{pmatrix}
x_1+o(x(z'))(z-z_1) \\[10pt]
y_1+o(1)_{x(z'') \to 0}(z-z_1) \\[10pt]
z
\end{pmatrix},
\end{equation}
\noindent where $z'$ and $z''$ are intermediate values between $z_1$ and $z$, and $h_{11}(\cdot)=:x(\cdot)$.
\par{}
We now show that we can replace $o(x(z'))$ by $o(x_1)$ in \eqref{eq:WsQ_14}. Let $z(s)=z_1+(z-z_1)s$ and $X(s)=x(z(s))=h_{11}(z(s))$, where $s \in [0,1]$. Our goal is to prove
\begin{equation*}
\dfrac{X(s)-x_1}{x_1}\to 0 \quad\mbox{as}\quad x_1 \to 0\quad \mbox{uniformly}.
\end{equation*}
Now suppose that there exist some $s_0$ and $\varepsilon>0$ such that $({X(s_0)-x_1})/{x_1} \geqslant \varepsilon$ for all $x_1$. In what follows we show a contradiction.
\par{} Note that we have $({X(0)-x_1})/{x_1}=0$. By the continuity of $X(s)$, there exists $s^* \in (0,s_0]$ such that $({X(s^*)-x_1})/{x_1} = \varepsilon$ and $({X(s)-x_1})/{x_1} < \varepsilon$ for every $s<s^*$. This further implies that, for any $s<s^*$, we have
\begin{equation}\label{eq:WsQ_n1}
X(s)=x_1+\varepsilon'(s) x_1,
\end{equation}
\noindent where $\varepsilon'(s)\in(0,\varepsilon)$ is a continuous function defined on $s\in[0,s^*]$. Equation \eqref{eq:WsQ_n1} along with (\ref{eq:WsQ_13}) leads to
\begin{equation*}
\dfrac{\mbox{d}h_{11}(z(s))}{\mbox{d}z}=o(h_{11}(z(s)))=o(X(s))=o(x_1+\varepsilon' x_1). 
\end{equation*}
\noindent Therefore, we have
\begin{equation*}
\dfrac{\;\;\dfrac{\mbox{d}h_{11}(z(s))}{\mbox{d}z}\;\;}{h_{11}(z(s))} =o(1)_{x_1\to 0},
\end{equation*}
\noindent and, particularly,
\begin{equation}\label{eq:WsQ_16}
\dfrac{\;\;\dfrac{\mbox{d}h_{11}(z(s))}{\mbox{d}z}\;\;}{h_{11}(z(s))} < \dfrac{\varepsilon}{2(z-z_1)s(1+\varepsilon')},
\end{equation}
\noindent for all $s\in(s_1,s^*)$ with any given $s_1\in(0,s^*)$ by choosing $x_1$ sufficiently small.
\par{} Let us now look at $X'(s)$ on $[0,s^*)$ given by
\begin{equation*}
X'(s)=\dfrac{\mbox{d}h_{11}(z(s))}{\mbox{d}z}(z-z_1)
\end{equation*}
\noindent which, by taking integral on both sides, yields 
\begin{equation}\label{eq:WsQ_15}
\begin{array}{rl}
X(s)-x_1=&(z-z_1) \displaystyle\int_0^s \dfrac{\mbox{d}h_{11}(z(s))}{\mbox{d}z} \mbox{d}s \\[20pt]
=&(z-z_1)s \dfrac{\mbox{d}h_{11}(z(s'))}{\mbox{d}z},
\end{array}
\end{equation}
\noindent where $s'\in (0,s)$. By plugging equation (\ref{eq:WsQ_15}) into (\ref{eq:WsQ_16}) and using \eqref{eq:WsQ_n1}, we obtain 
\begin{equation*}
\dfrac{\;\;\dfrac{X(s)-x_1}{(z-z_1)s}\;\;}{x_1+\varepsilon' x_1} < \dfrac{\varepsilon}{2(z-z_1)s(1+\varepsilon')},
\end{equation*}
\noindent i.e.
\begin{equation}\label{eq:WsQ_17}
X(s)<x_1+\dfrac{\varepsilon}{2} x_1,
\end{equation}
\noindent which holds for all $s\in(s_1,s^*)$. This contradicts the continuity of $X(s)$ since we have $\lim_{s\to s^*}X(s)=x_1+\varepsilon x_1$.
\par{}
We now have proved $o(x(z'))\sim o(x_1)$. By a similar argument, we also have $o(x(z''))\sim o(x_1)$. The function of $W^{ss}(M_1)$ arrives at the following form:
\begin{equation}\label{eq:WsQ_18}
h_1(z)=
\begin{pmatrix}
h_{11}(z)\\
h_{12}(z)\\
z
\end{pmatrix}
=
\begin{pmatrix}
x_1+o(x_1)(z-z_1) \\[10pt]
y_1+o(1)_{x_1 \to 0}(z-z_1) \\[10pt]
z
\end{pmatrix}.
\end{equation} 
\noindent which leads to the statement of Lemma \ref{lem:WsQ}.\qed
\section{Proof of Theorem \ref{thm:hetero_2}}\label{sec:prf2}
We first prove that, by an arbitrarily small perturbation to system $X$, we can simultaneously find heterodimensional cycles and a wild hyperbolic set. Then, we show that with condition (C5) satisfied the heterodimensional cycles and the wild hyperbolic set belong to the attractor $\mathcal{A}$.
\begin{lem}\label{lem:hdcwhs}
There exists a sequence $\{X_m\}$ of systems converging to system $X$ in $C^r$ topology such that each system $X_{m}$ contains heterodimensional cycles as well as a wild hyperbolic set.
\end{lem}
\noindent{\it Proof.} By Shilnikov theorem (see Lemma \ref{lem:WsP}), there exists a countable set of periodic orbits of $X$ with index-1 in any small neighbourhood of the homoclinic loop $\Gamma^+$. Besides, it is shown in \citep{sh70} that, for some sufficiently large integer $K>0$ and any $\rho'$ close to $\rho^*$, there exists an invariant hyperbolic set $\Lambda_{K,\rho'}$ in any such neighbourhood with one-to-one correspondence to the set of two-sided sequences $\{i_n\}_{-\infty}^{+\infty}$, where $\rho'i_n\leqslant i_{n+1}$ and $i_n\geqslant K$ for all $n$. For any small $\mu\neq 0$, there exists $\bar{K}\gg K$ such that one can find a closed invariant hyperbolic set $\Lambda_{K,\bar{K},\rho'}$ with one-to-one correspondence to the set of two-sided sequence $\{i_n\}_{-\infty}^{+\infty}$, where $\rho'i_n\leqslant i_{n+1}$ and $K\leqslant i_n\leqslant \bar{K}$ for all $n$. The purpose of finding $\bar{K}$ is to single out a closed subset $\Sigma_{K,\bar{K},\rho'}$ from $\Sigma_{K,\rho'}$ such that $\Lambda_{K,\bar{K},\rho'}$ survives from small perturbation due to its closeness. Note that the sets $\Lambda_{K,\rho'}$ and $\Lambda_{K,\bar{K},\rho'}$ corresponding to different $\rho'$ are different, but they all exist in system $X$. We drop the subscript $\rho'$ of these sets for avoiding ambiguity.
\par{}
The set $\Lambda_{K,\bar{K}}$ can be wild. Indeed, Theorem 1 of \citep{os87} states that there exists a dense set $\{\bar{\rho}_m\}$ of $\rho$ values in $(0,{1}/{2})$ such that system $X_{\bar{\rho}_m}$ contains a homoclinic tangency associated to a periodic orbit $L$ in $\Lambda_{K}$, and the original homoclinic loops are kept. By choosing $\bar{K}$ sufficiently large, we can ensure that the periodic orbit $L$ lie in the set $\Lambda_{K,\bar{K}}$, and therefore make $\Lambda_{K,\bar{K}}$ a wild hyperbolic set.
\par{}
By Newhouse theorem (see \citep{gts93,pv94}), for each $\bar{\rho}_m$, we have a small neighbourhood of it such that, for each $\rho$ value in it, the corresponding systems $X_{\rho}$ contains a wild hyperbolic set. Moreover, systems sufficiently close to system $X_\rho$ in $C^r$ topology also contain such sets. Let $\varepsilon$ be the set of parameters other than $\rho$. Let $\rho*$ be the $\rho$ value of the original system $X$. It follows that there exist infinitely many open neighbourhood $B_{\bar{\rho}_m}$ of $X_{\bar{\rho}_m}$ in the space of $R$-symmetric systems where $|\bar{\rho}_m-\rho|<1/2^m$ such that any system in these neighbourhoods has a wild hyperbolic set $\Lambda_{K,\bar{K}}$. 
\par{}
Now we can apply Theorem \ref{thm:hetero_1} to system $X_{\bar{\rho}_m}$ since it satisfies all conditions required by the theorem. Thus, in each ball $B_{\bar{\rho}_m}$, we can find a system $X_m$ such that it  contains heterodimensional cycles as well as a wild hyperbolic set. \qed    

%
%
%
%
%
\par{}
We next prove that the set $\Lambda_{K,\bar{K}}$ and the heterodimensional cycle coexist in the attractor $\mathcal{A}$ of system $X_m$ with sufficiently large $m$. Recall the cross-section defined in Section \ref{sec:2.2}, which is
\begin{equation*}
S=\{(x_1,x_2,y,z)|\|(x_1,x_2)\|=1,|y|\leqslant 1,\|z\|\leqslant 1\}.
\end{equation*}
\noindent It is also a cross-section for $W^s_{loc}(O)$. Denote $W^s_{loc}(O)\cap S$ by $S_0$, $\{y\geqslant 0\}\cap S$ by $S_1$ and $\{y\leqslant 0\}\cap S$ by $S_2$. By the assumption of the volume hyperbolicity in $D$, system $X$ always have a strange attractor $\mathcal{A}$ in $D$. 
\par{}
Let us recall some properties of this attractor given by Theorem 3 in \citep{ts98}. Denote by $M^{\pm}_i$ the successive intersections of separatrices $\Gamma^\pm$ with the cross-section $S$. The intersection $\mathcal{A}\cap S$ consists of $N$ connected components, where $2\leqslant N <2+l(\rho)$ is a finite number bounded above by a function $l(\rho)>0$ with $0<\rho<{1}/{2}$, and each connected component contains at least one of the points $M^{\pm}_i$. Note that $N$ is independent of parameters other than $\rho$. Moreover, we have two integers $N^+$ and $N^-$ satisfying $N^+ +N^-=N$ such that 
\begin{equation*}
\mathcal{A}\cap S=A_1^+\cup\dots\cup A_{N^+}^+\cup A_1^-\cup\dots\cup A_{N^-}^-,
\end{equation*}
\noindent where $A_i^+$ and $A_j^-$ are disjoint connected components. Denote by $T_S\equiv(T_{S_1},T_{S_2})$ the Poincaré map on $S=S_1\cup S_2$. We have 
\begin{equation}\label{eq:connected_components_1}
A^+_i\cap S_0=\emptyset ,\quad A^-_j \cap S_0=\emptyset
\end{equation}  
\noindent and
\begin{equation}\label{eq:connected_components_2}
\begin{array}{c}
T_{S_1}((A^+_{N^+}\cup A^-_{N^-})\cap S_1)=A^+_1 ,\quad\quad
T_{S_2}((A^+_{N^+}\cup A^-_{N^-})\cap S_2)=A^-_1 , \\[10pt]
A^+_i=T_{S_1}A^+_1, \quad A^-_j=T_{S_2}A^-_1 ,
\end{array}
\end{equation}
\noindent where $1<i<N^+$ and $1<j<N^-$.
\par{}
We assume that, at $\mu=0$, the unstable separatrices $\Gamma^+$ and $\Gamma^-$ of the system $X$ intersect $S$ for $M$ times in total such that $M>N$, where $N$ is the above-mentioned number of the connected components. This assumption is not a restriction. Indeed, results in \citep{fe93} state that a homoclinic loop to a saddle focus equilibrium can be split in a way such that the splitting forms a new loop with arbitrary more rounds with respect to the original one. Suppose now the system $X$ has a single round homoclinic loop. Let us choose a surface $L$ corresponding to systems having a homoclinic loop to a saddle-focus equilibrium in the space of dynamical systems, and let it contain the system $X$. Then the surfaces corresponding to systems having homoclinic loops with arbitrarily many rounds are accumulating on $L$. Since our family $X_\varepsilon$ is transverse to $L$ with respect to $\mu$, it is transverse to all close surfaces. Let $X_{\varepsilon_0}:=X$. It follows that we can pick $X_{\varepsilon^*}$, with $\varepsilon^*$ arbitrarily close to $\varepsilon_0$ from either left or right, such that the system at $\varepsilon=\varepsilon^*$ has a homoclinic loop with more than one round, which intersects $S$ for multiple times. Any result on bifurcation that holds for all such system $X_{\varepsilon^*}$ will automatically hold for $X_{\varepsilon_0}$. 
\par{}
We now choose the cross-section $\Pi$ used in the proof of Theorem \ref{thm:hetero_1} sufficiently small such that the loops $\Gamma^\pm$ are still single-round with respect to $\Pi$ in system $X$, and therefore double-round in systems $X_m$. This makes all the arguments used for proving Theorem \ref{thm:hetero_1} valid here. 
\par{}
Let $P_1$ be the index-1 point of the Poincaré map $T$ on $\Pi$ used to create a heterodimensional cycle in system $X_m$ by Theorem \ref{thm:hetero_1}. In system $X_m$, the periodic orbit $L$ of the flow whose stable manifold intersects its unstable manifold non-transversely is double-round with respect to $\Pi$. There are two points in $L\cap\Pi$, and we denote by $P_2\in\Pi_1$ the one closer to $\Pi_0=\Pi\cap\{x=0\}$. By Lemma \ref{lem:ubdyk*}, the point $P_1$ can be chosen such that its $x$-coordinate $x_1$ is in $(0,|\mu|)$. The point $P_2$ can be arbitrarily close to $\Pi_0$ by choosing $L$ close to $\Gamma^+$ (see the proof of Theorem 1 in \citep{os87}). Especially, we can let its $x$-coordinate $x_2$ also lie in $(0,|\mu|)$. By Lemma \ref{lem:WsP}, we know that the local stable manifolds $W^s_{loc}(P_i)$ $(i=1,2)$ are given by
\begin{equation}
x=x_i+o(1)_{x_i\to 0},
\end{equation}
\noindent which are bounded between $\Pi_0$ and $\{x=|\mu|\}$. We now follow the backward orbits of the flow starting on $W^s_{loc}(P_i)$ until they intersect the cross-section $S$. This gives us two periodic points $P'_i\in S$ of the map $T_S$ along with their local stable manifolds $W^s_{loc}(P_i')$. One can check that those manifolds are bounded by $S_0$ and $\{x=o(1)_{\mu\to 0}\}$, where the small term is positive and is a continuous function of $x_1,x_2$ and $\mu$. In fact, it is shown in \citep{ts98} that the stable manifold of a fixed point of $T_S$ is a nearly horizontal curve circling around $O$ on the cylinder-like cross-section $S$; by choosing a small part $S^*$ from $S$ such that $P'_i$ is a fixed point of $T_S|_{S^*}$, we can achieve the same result on $W^s(P'_i)$. This means that points on $W^s(P'_i)$ have the same $x$-coordinate $x'_i$ of $P'_i$ up to some small corrections. Obviously, $x'_i\to 0$ as $x_i\to 0$. Note that we can pick $P_i$ with sufficiently small $x$-coordinates when $|\mu|$ is sufficiently small. Therefore, the stable manifolds $W^s(P'_i)$ $(i=1,2)$ lie under a surface $\{x=o(1)_{\mu\to 0}\}$.
\par{} 
We are now in the position to finish the proof. Recall the discussion on the attractor $\mathcal{A}$. Since we have $M>N$, the connected component $M^+_{N^+}\in A^+_{N^+}$ is at a finite distance from $S_0$. Therefore, we obtain a connected curve $l^+\subset A^+_{N^+}$ joining $M^+_{N^+}$ and a point on $S_0$. If $l^+\subset S_1$, then we have $W^s_{loc}(P_i')\cap l\neq\emptyset$ by choosing $\mu$ sufficiently small; if $l^+\subset S_2$, then we achieve the same result by consider similar points $P_i$ from $\Pi_2$. Note that the manifolds $W^s_{loc}(P_i')$ are sufficiently long to intersect the connected components since they circle around $O$ on $S$ for multiple times. Recall that attractor $\mathcal{A}$ is the set of points accessible from the equilibrium $O$. Hence, the set $\mathcal{A}\cap S$ contains $W^s(P_i')$ and $P_i'$ along with their unstable manifolds, which implies that it contains a heterodimensional cycle as well as a wild hyperbolic set of the map $T$. When we return to the full system $X_m$, we find a heterodimensional cycle and the wild hyperbolic set $\Lambda_{K,\bar{K}}$ coexisting in $\mathcal{A}$. The theorem is proved.

\begin{thebibliography}{1}
\bibitem[An67]{an67} Anosov, D. V.. Geodesic flows on closed Riemannian manifolds of negative curvature.
Proc. Steklov Inst. Math., 90 (1967), 1–235.
%
\bibitem[ABS77]{abs77} Afraimovich, V. S., Bykov, V. V. and Shilnikov, L. P.. On the origin and structure of the Lorenz attractor. Akademiia Nauk SSSR Doklady, 234 (1977), 336-339.
%
\bibitem[ABS83]{abs83} Afraimovich, V. S., Bykov, V. V. and Shilnikov, L. P.. On the structurally unstable attracting limit sets of Lorenz attractor type. Tran. Moscow Math. Soc., 2 (1983), 153-215. 
%
\bibitem[BC15]{bc15} Bonatti, C. and Crovisier, S.. Center manifolds for partially hyperbolic sets without strong unstable connections. Journal of the Institute of Mathematics of Jussieu, available on CJO2015. doi:10.1017/S1474748015000055. 
%
\bibitem[BD96]{bd96} Bonatti, C. and Díaz, L. J.. Persistent transitive diffeomorphisms. Annals of Mathematics, 143(2) (1996), 357-396.
\bibitem[BDV00]{bdv} Bonatti, C., Díaz, L. J. and Viana, M.. Dynamics Beyond Uniform Hyperbolicity. Springer, Berlin, Heidelberg, New York, 2000.
\bibitem[BD08]{bd08} Bonatti, C. and Díaz, L. J.. Robust heterodimensional cycles and C1-generic dynamics. J. Inst. Math. Jussieu 7, no. 3 (2008), 469-525.
\bibitem[BSS12]{bss12} Barrio, R., Shilnikov, A. L. and Shilnikov L. P.. Kneadings, Symbolic Dynamics and Painting Lorenz Chaos. International Journal of Bifurcation and Chaos, 22(4) (2012), 1230016(1) - 1230016(24). 
\bibitem[Dí92]{d92} Díaz, L. J. and Rocha, J.. Non-connected heterodimensional cycles: bifurcation and stability. Nonlinearity, 5 (1992), 1315-1341.
\bibitem[Dí95a]{d95} Díaz, L. J.. Robust nonhyperbolic dynamics and heterodimensional cycles. Ergodic Theory and Dynamical Systems, 15 (1995), 291-315.
\bibitem[Dí95b]{d95_2} Díaz, L. J.. Persistence of cycles and nonhyperbolic dynamics at the unfolding of heteroclinic bifurcations. Nonlinearity, 8 (1995), 693-715.
\bibitem[EFF82]{eff82} Evans, J. W., Fenichel, N. and Feroe, J. A.. Double impulse solutions in nerve axon equations. SIAM J. Appl. Math. 42 (1982) 219.
\bibitem[Fe93]{fe93} Feroe, J. A.. Homoclinic orbits in a parametrized saddle-focus system, Phys. D 62:1-4 (1993), 254-262.
\bibitem[Ga83]{ga83} Gaspard, P.. Generation of a countable set of homoclinic flows through bifurcation, Physics Letters A, Volume 97, Issues 1-2, 8 August 1983, 1-4.
\bibitem[GST93]{gts93} Gonchenko, S. V., Turaev, D. V. and Shilnikov, L. P.. On the existence of Newhouse regions in a neighborhood of systems with a structurally unstable homoclinic Poincaré curve (the multidimensional case). (Russian) Dokl. Akad. Nauk 329 (1993), no. 4, 404--407; translation in Russian Acad. Sci. Dokl. Math. 47 (1993), no. 2, 268–273.
%
\bibitem[GST09]{} Gonchenko, S. V., Shilnikov, L. P. and Turaev, D. V.. On global bifurcations in three-dimensional diffeomorphisms leading to wild Lorenz-like attractors. Reg Chaot Dyn, 14:1 (2009), 137.
%
\bibitem[GTGN97]{gtgn97} Gonchenko, S. V., Turaev, D. V., Gaspard, P. and Nicolis, G.. Complexity in the bifurcation structure of homoclinic loops to a saddle-focus, Nonlinearity, 10 (1997), 409-423.
%
\bibitem[HPS77]{hps77} Hirsch, M., Pugh, C. and Shub, M. Invariant manifolds. Springer-Lecture Notes on Mathematics, 583, Heidelberg, 1977.
%
\bibitem[HS10]{hs10} Homburg, A. J. and Sandstede, B.. Homoclinic and heteroclinic bifurcations. Vector Fields. Handbook of Dynamical Systems III, Elsevier, (2010), 379-524.
%
\bibitem[Hu82]{hu82} Hurley, M.. Attractors: persistence and density of their basins. Trans. Amer. Math. Soc., 269 (1982), 247-271.
%
\bibitem[Li16]{li16} Li, D.. Homoclinic bifurcations that give rise to heterodimensional cycles near a Saddle-focus equilibrium. arXiv:1604.00431 [math.DS].
%
\bibitem[Ne79]{newhouse1} Newhouse, S. E.. The abundance of wild hyperbolic sets and non-smooth stable sets for diffeomorphisms. Inst. Hautes Études Sci. Publ. Math, 50 (1979), 101-151.
%
\bibitem[NP76]{np76} Newhouse, S. and Palis, J.. Cycles and bifurcation theory. Asterisque, 31 (1976), 43–140.
%
\bibitem[OS87]{os87} Ovsyannikov, I. M. and Shilnikov, L. P.. On systems with a saddle-focus homoclinic curve. Sbornik: Mathematics, 58(2) (1987), 557-574.
%
\bibitem[OS92]{os92} Ovsyannikov, I. M. and Shilnikov, L. P.. Systems with a homoclinic curve of multidimensional saddle-focus type, and spiral chaos. Math. USSR Sbornik, 73 (1992),  415-443.
%
\bibitem[PV94]{pv94} Palis, J. and Viana, M.. High dimension diffeomorphisms displaying infinitely many periodic attractors. Ann. of Math., (2) 140 (1994), no. 1, 207–250.
%
\bibitem[Pa00]{p00} Palis, J. A global view of dynamics and a conjecture on the denseness of finitude of attractors. Géométrie complexe et systèmes dynamiques (Orsay, 1995). Astérisque No. 261 (2000), xiii–xiv, 335–347.

\bibitem[Ru81]{ru81} Ruelle, D.. Small random perturbations of dynamical systems and the definition of attractors. Comm. Math. Phys., 82 (1981), 137-151. 
%
\bibitem[ST99]{st99} Shashkov, M. V. and Turaev, D. V.. An Existence theorem of smooth nonlocal center
manifolds for systems close to a system with a
homoclinic loop. J. Nonlinear Sci., 9 (1999), 525-573.
%
\bibitem[Sh65]{sh65} Shilnikov, L. P.. A case of the existence of a countable number of periodic motions (Point mapping proof of existence theorem showing neighborhood of trajectory which departs from and returns to saddle-point focus contains denumerable set of periodic motions). SOVIET MATHEMATICS, 6 (1965), 163-166.
%
\bibitem[Sh70]{sh70} Shilnikov, L. P.. A contribution to the problem of the structure of an extended neighborhood of a rough equilibrium state of saddle-focus type. Sbornik: Mathematics, 10(1) (1970), 91-102.
%
\bibitem[SSTC01]{sstc1} Shilnikov, L. P., Shilnikov, A. L., Turaev, D. V. and Chua, L. O.. Methods Of Qualitative Theory In Nonlinear Dynamics (Part I). World Sci.-Singapore, New Jersey, London, Hong Kong, 2001.
%
\bibitem[SSTC01]{sstc2} Shilnikov, L. P., Shilnikov, A. L., Turaev, D. V. and Chua, L. O.. Methods Of Qualitative Theory In Nonlinear Dynamics (Part II). World Sci.-Singapore, New Jersey, London, Hong Kong, 2001.
%
\bibitem[Tu02]{tu02} Tucker, W.. A rigorous ODE solver and Smale’s 14th problem. Foundations of Computational Mathematics, 2(1) (2002), 53-117.
%
\bibitem[Tu96]{tu96} Turaev, D. V.. On dimension of non-local bifurcational problems, International Journal of Bifurcation and Chaos, 6(5) (1996), 919-948.
%
\bibitem[TS98]{ts98} Turaev, D. V. and Shilnikov, L. P.. An example of a wild strange attractor. Sbornik. Math. 189(2) (1998), 291-314.
\end{thebibliography}
\end{document}